\documentclass[11pt,a4paper]{article}
\usepackage{etex}
\usepackage{inputenc}
\usepackage[english]{babel}
\usepackage{lmodern}					% harmonie des polices
\usepackage{textcomp}				% table de symboles
\usepackage{verbatim}			% pour faire des commentaires par block avec begin{comment}...end{comment}
\usepackage[toc,page]{appendix} 		% pour faire des annexes
\usepackage{cite}
\usepackage{amsmath,					% packages mathmatiques
			amsthm,
			amsfonts,
			amssymb, 
			mathtools}%,a4wide}
\usepackage{stmaryrd}
\usepackage{physics}
\usepackage{enumitem}
\usepackage{hyperref}
\hypersetup{
colorlinks=true,
linkcolor=purple,
filecolor=magenta,
urlcolor=cyan,
citecolor=blue
}

 \usepackage{tikz-cd} 
\usetikzlibrary{tikzmark}
\usetikzlibrary{arrows}
 \usepackage{verbatim}
\definecolor{DarkGreen}{HTML}{1cad22}
\usepackage[hmargin=3cm,vmargin=2.5cm]{geometry}
\usepackage{graphicx} % Required for inserting images
\usepackage{amsthm, amssymb, amsmath}
\usepackage{color}
\usepackage{xcolor}
\usepackage{tikz}
\usepackage{tikz-cd}
\usepackage{amsmath}
\usepackage{amssymb}
\usepackage{float}
\usepackage{caption}
\usepackage{mathtools}
\usepackage{hyperref}
\usepackage{mathrsfs}

\usepackage{graphicx}
\usepackage{epsfig}
\newcommand{\HH}{\ensuremath{\mathbb H}}

\newcommand{\ZZ}{\ensuremath{\mathbb Z}}

\numberwithin{equation}{section}

{\theoremstyle{definition}\newtheorem{definition}{Definition}[section]

\newtheorem{defititle}[definition]{\Title}

\newtheorem{remark}[definition]{Remark}
\newtheorem{example}[definition]{Example}

}
\newtheorem{prop}[definition]{Proposition}
\newtheorem{proposition-definition}[definition]{Proposition-Definition}
\newtheorem{lemma}[definition]{Lemma}
\newtheorem{thm}[definition]{Theorem}
\newtheorem{cor}[definition]{Corollary}

\newtheorem*{theorem}{Theorem}

\newtheoremstyle{named}{}{}{\itshape}{}{\bfseries}{.}{.5em}{\thmnote{#3's }#1}
\theoremstyle{named}

\usepackage{amsthm}

\newenvironment{mainthm-env}[1]
  {\begin{theorem}[Main Theorem #1]}
  {\end{theorem}}

\title{The topology of local quaternionic toric actions}

\author{Panagiotis Batakidis$^\dagger$ \\\href{mailto:batakidis@math.auth.gr}{batakidis@math.auth.gr}  \and Ioannis Gkeneralis$^\dagger$
\\\href{mailto:igkeneralis@math.auth.gr}{igkeneralis@math.auth.gr}}

\date{$^\dagger$ Department of Mathematics, Aristotle University of Thessaloniki}

\newcommand{\keywords}[1]{\textbf{Keywords: } #1}
\newcommand{\msc}[1]{\textbf{Mathematics Subject Classification: } #1}

\begin{document}

\maketitle

\begin{abstract}
In this paper we examine the topology of manifolds equipped with a local quaternionic toric action modeled on the regular representation of the quaternionic torus $Q^n=(S^3)^n$. Building on our previous work, where the toric, differential and tetraplectic foundations were established, we show that the global topology of such manifolds is determined by the orbit space and its characteristic data. We construct Leray--Serre and Atiyah--Hirzebruch spectral sequences for the orbit projection, yielding explicit descriptions of the cohomology and $K$-theory of manifolds equipped with local quaternionic toric actions. In dimension four, we develop a quaternionic analogue of the Meyer signature formula and we briefly outline an $L$-theoretic interpretation of the resulting signature invariants. These results extend the methods of the classical (complex) toric topology to the quaternionic setting.
\end{abstract}

% Keywords and MSC
\noindent \keywords{local quaternionic toric action, regular representation, fundamental group, cohomology, Leray/Serre spectral sequence, untwisted $K$-theory, signature} \\[1ex]
\msc[2020] {57S25, 55T10, 57R18, 57N65}

\vspace{0.5ex}
\noindent \textbf{Conflict of Interest and Data Availability Statement:} The authors state that there is no conflict of interest to declare. Moreover, data sharing is not applicable to this article as no datasets were generated or analyzed during the research carried out in this paper.

\tableofcontents

\section{Introduction}

In toric topology, one studies the topology of manifolds with torus actions using the combinatorics of their orbit spaces. The typical example, due to Davis and Januszkiewicz \cite{dj}, is that of \emph{quasitoric manifolds}: smooth $2n$--manifolds endowed with a locally standard $(S^1)^n$--action whose orbit (quotient) space is a simple convex polytope. In this framework, the global topology of the manifold is completely determined by combinatorial and lattice data attached to the facets of the polytope, giving rise to explicit algebraic models for cohomological and $K$--theorical descriptions (for details see \cite{BuchPan2000,bp,tt}). These ideas have since shaped a large part of modern toric and transformation group topology.

\medskip

A natural problem is to understand how this correspondence extends to  the non-commutative setting. The quaternionic analogue of this picture arises when one replaces the complex torus $(S^1)^n$ by the compact quaternionic torus,
\[
Q^n=(S^3)^n\simeq SU(2)^n \simeq \mathrm{Sp}(1)^n,
\]
acting on $\mathbb{H}^n$ by coordinatewise left quaternionic multiplication. This representation provides the basic noncommutative analogue of the complex toric model and is called the \emph{regular representation}. Manifolds locally modeled on this action—appearing in the works of Scott~\cite{scott1995} and Hopkinson~\cite{ho} under the name \emph{quaternionic toric varietes} and \emph{quoric manifolds} respectively—share much of the formal structure of quasitoric manifolds, but with new phenomena arising from the noncommutativity of quaternions.

\medskip

In previous work \cite{bg25}, we introduced the notion of a \emph{local quaternionic toric action} modeled on the regular representation of $(S^3)^n$ and studied its compatibility with tetraplectic structures. There, we developed the analytic and differential-geometric foundations of the theory, including the topological classification of local models, the construction of canonical orbit spaces with corners, and the description of the corresponding characteristic bundles which encode the local linear (isotropy) data. This provided a quaternionic extension of the ideas introduced by Yoshida's theory \cite{Yoshida} of locally standard torus actions.

\medskip

The purpose of the present paper is to investigate the \emph{global topological invariants} of manifolds endowed with local quaternionic toric actions, building on the geometric framework established in \cite{bg25}. Our main objective is to show that the machinery of spectral sequences used in toric topology extends naturally to the quaternionic setting and provides concrete computational tools for fundamental groups, cohomology, $K$-theory, and signature-type invariants.

\medskip

More precisely, let a compact connected  \(4n\)-dimensional manifold $M^{4n}$ equipped with a local $Q^n$-action and orbit map
\[
\pi_M : M^{4n} \longrightarrow B_M.
\]
As shown in \cite{bg25}, the orbit space $B_M$ is an $n$-dimensional manifold with corners endowed with a natural quaternionic integral affine structure. We show that much of the global topology of $M^{4n}$ is encoded in the orbit space $B_M$. The transition automorphisms $\{\rho_{\alpha\beta}\}\subset \mathrm{Aut}(Q^n)$ of a weakly regular atlas (Definition \ref{def:weakly_regular}) determine a principal $\mathrm{Aut}(Q^n)$-bundle $P_M\to B_M$ and an associated lattice subbundle $\mathcal{L}_M$ over the codimension-one strata of $B_M$. The pair $(P_M,\mathcal{L}_M)$, is called the \emph{characteristic pair} (Definition \ref{def: characteristic_pair}), from which we construct the associated \emph{canonical model}  $M_{(P_M,\mathcal{L}_M)}$ (Definition \ref{def:canonical_model}) that provides a topological representative of the local action.

\medskip

Building on this local model of quaternionic toric actions and the associated stratified orbit projection \(\pi_M:M^{4n}\to B_M\), we investigate the fundamental topological invariants of \(M^{4n}\). From this geometric framework, we derive a series of structural results connecting the topology of the total space to that of its orbit space.

The orbit projection $\pi_M$ induces an isomorphism on fundamental groups, reflecting the fact that the quaternionic torus fibers are simply connected. More precisely, in Proposition \ref{prop:pi1}, we prove that for a compact connected manifold \(M^{4n}\) with a local \(Q^n\)-action, the orbit map induces an isomorphism
\[
\pi_1(M^{4n}) \;\cong\; \pi_1(B_M)
\]
and then we provide a construction to compute cohomology groups. The Leray--Serre type spectral sequence associated with the orbit filtration
\[
M^{(p)}=\pi_M^{-1}(B_M^{(p)})
\]
provides a natural tool for computing \(H^\bullet(M^{4n};\mathbb Z)\). Its cellular form takes into account the stratified nature of the orbit projection: over a cell contained in a face \(F\subset B_M\), the relevant coefficient group is the cohomology of the corresponding orbit \(Q^F\). In low degrees, the simply connectedness of the quaternionic orbits implies the expected identifications
\[
H^1(M^{4n};\mathbb Z)\cong H^1(B_M;\mathbb Z),
\qquad
H^2(M^{4n};\mathbb Z)\cong H^2(B_M;\mathbb Z)
\]
under the untwisted hypotheses considered below. Therefore, we have:

\begin{theorem}[Theorems \ref{thm:leray} \& \ref{thm:cellularSSS}]
There exists a first–quadrant spectral sequence \((E_M)^{p,q}_r\) converging to the cohomology of $M^{4n}$, \(H^{p+q}(M^{4n};\mathbb{Z})\), in the sense that
\[ E_\infty^{p,q} \;\cong\; F^pH^{p+q}(M^{4n};\mathbb{Z})/F^{p+1}H^{p+q}(M^{4n};\mathbb{Z}),\] for the natural filtration $F^\bullet$ induced by the skeleta of $B_M$.
\end{theorem}

As an illustration (in paragraph \ref{par:quoric}) we apply the cellular spectral sequence to quoric manifolds \cite{ho}. For \(M=\mathbb H P^1\cong S^4\), with orbit space
\(\Delta^1\), the calculation recovers
\[
H^\ast(\mathbb H P^1;\mathbb Z)\cong \mathbb Z[u]/(u^2),
\qquad |u|=4.
\]
More generally, our computation for the examples \(M=\mathbb H P^n\) is compatible with Hopkinson's formula (see \cite[\S 6]{ho}), that
\[
H^\ast(\mathbb H P^n;\mathbb Z)\cong \mathbb Z[u]/(u^{n+1}),
\qquad |u|=4.
\]

\medskip

An analogous statement holds in complex \(K\)-theory: replacing ordinary cohomology by \(K^\bullet\) in the same orbit-filtration construction gives an Atiyah-Hirzebruch type spectral sequence computing \(K^\bullet(M)\). Its first page is given by
\[
(E_M)^{p,q}_1
=
K^{p+q}(M^{(p)},M^{(p-1)}),
\]
and admits a cellular description in terms of the \(K\)-groups of the orbit fibers over the cells of \(B_M\); see Theorem~\ref{thm:K-theory-SS}.

In real dimension four, we also extend Novikov additivity to the quaternionic/ noncommutative setting and define a \emph{quaternionic Meyer cocycle} (Theorem \ref{thm:quat_meyer}) that governs the signature of $M$. This provides a quaternionic counterpart to the signature formula for toric $4$--manifolds and links the algebraic data of the local action to the intersection form of the manifold. 

The appearance of a quaternionic Meyer cocycle reflects a genuinely noncommutative feature of the theory. In contrast with the complex toric case, where the torus $(S^1)^n$ is abelian and the corresponding monodromy is encoded purely by lattice data, the quaternionic torus $Q^n$ is built from the nonabelian group $Q$. As a result, the local isotropy representations and their monodromies carry additional information which influences the intersection form. The quaternionic Meyer cocycle detects precisely this extra twisting phenomenon, which is invisible in purely combinatorial models.

Beyond the cohomological and $K$--theoretic structure, we also outline an $L$--theoretic interpretation of the quaternionic signature. From this viewpoint, the signature of a local quaternionic toric $4$--manifold represents the degree--$4$ component of its $L$--class in $L^\bullet(\mathbb{Z})$, while potential nontrivial Euler classes or monodromies of the quaternionic torus bundle correspond to higher signatures in the sense of Novikov. This observation places the quaternionic toric framework within the broader landscape of topological invariants derived from generalized cohomology theories, where $K$--theory captures the multiplicative, representation--theoretic features of the fiber, and $L$--theory encodes the additive, intersection--theoretic aspects of the total space.

\medskip

Together, these results establish the topological foundation of local quaternionic toric actions and demonstrate that their global invariants are determined by the orbit space and the characteristic data $(P_M,\mathcal{L}_M)$. 

\medskip
Our results apply in particular to globally defined quaternionic toric manifolds in the sense of Gentilli–Gori–Sarfatti \cite{ggs}, to quoric manifolds of Hopkinson \cite{ho}, and to Scott’s singular models \cite{scott1995}. In each case the action is locally modeled on the regular representation of $Q^n$ (see Subsection \ref{sec:loc_quat_toric_mfds}), so that these examples fit naturally into the present framework.

\medskip

 The approach developed here extends the methods of \cite{Yoshida} from the classical (complex) toric topology to the quaternionic toric topology. The main difficulty arises when dealing with noncommutativity of the group actions, spectral sequences, and bundle-theoretic structures. As in Yoshida’s work, we focus primarily on the global case in which the Euler class vanishes, so that the orbit map admits a section and the associated spectral sequences may be described in an untwisted form. This already includes all globally defined quaternionic toric manifolds known in the literature. In addition, however, we include remarks on the case of nonvanishing Euler class, where the resulting monodromy gives rise to genuinely twisted versions of these spectral sequences--a phenomenon that does not appear in Yoshida’s treatment.

\subsection*{Organization of the paper}

Section~\ref{sec:prelims} reviews the basic constructions and notation for local quaternionic toric actions, including regular atlases, orbit spaces, and characteristic pairs $(P_M,\mathcal{L}_M)$. We also introduce the canonical model associated with $(P_M,\mathcal{L}_M)$ and summarize its key topological properties. 

Section~\ref{sec:topology1} develops the topological analysis: we first determine the fundamental group of $M^{4n}$, showing that $\pi_1(M)\cong\pi_1(B_M)$, and then construct the Leray--Serre spectral sequence computing the cohomology of $M$.

In Section~\ref{sec:topology2} we construct the corresponding cellular spectral sequence in complex \(K\)-theory, obtained by applying \(K^\bullet\) to the filtered pairs \((M^{(p)},M^{(p-1)})\). This gives the \(K\)-theoretic analogue of the cohomological construction. Finally, we focus on the four-dimensional case, where Novikov additivity and a quaternionic version of the Meyer cocycle yield a formula for the signature. We conclude with a brief discussion of the $L$--theoretic interpretation of the quaternionic signature, placing these results in the broader context of higher signatures and the Novikov conjecture.

\subsection*{Basic assumption}

Set \(e(M)=e(M,Q^n)\) to be the Euler class of \(M\), as described in \cite[\S 4.1]{bg25}. Throughout this paper we assume that \(e(M)\) vanishes unless otherwise stated. This is the untwisted case in the sense that the obstruction to the existence of a global section of the orbit projection vanishes. The orbit projection remains stratified, and the spectral sequences used below are formulated with respect to the orbit-type filtration
\[
M^{(p)}=\pi_M^{-1}(B_M^{(p)}).
\]
We briefly indicate how nonvanishing Euler classes would lead to twisted versions of these constructions - a topic that will be developed in future work.

\section{Preliminaries}
\label{sec:prelims}

This section collects the material needed to define \emph{local quaternionic toric manifolds} and to construct the associated canonical space.  We fix notation and state the constructions and key properties; proofs are omitted since they appear in detail in \cite{bg25}.

\subsection{Local quaternionic toric manifolds}\label{sec:loc_quat_toric_mfds}
Set \(Q=S^3\cong SU(2) \cong \mathrm{Sp}(1)\) and \(Q^n=(S^3)^n\) as the \emph{quaternionic torus}. We refer to the \emph{regular representation} of $Q^n$, as the regular action of $Q^n$ on $\mathbb{H}^n$, i.e. the componentwise quaternionic left-multiplication of $Q^n$ on $\mathbb{H}^n$.

\medskip
\begin{definition}[Regular coordinate neighborhood]
Let $M^{4n}$ be a 4n-dimensional smooth manifold equipped with a smooth $Q^n$-action. A \emph{regular coordinate neighborhood} of $M^{4n}$ is a triple $(U,\rho,\varphi)$ consisting of a $Q^n$-invariant open set $U$ of $M^{4n}$, an automorphism $\rho$ of $Q^n$, and a $\rho$-equivariant diffeomorphism $\varphi$ from $U$ to some $Q^n$-invariant open subset of $\mathbb{H}^n$. The action of $Q^n$ on $M^{4n}$ is said to be \emph{locally regular} if every point in $M^{4n}$ lies in some regular coordinate neighborhood, and an atlas consisting of regular coordinate neighborhoods is called \emph{regular atlas}.
\end{definition}

\medskip
\begin{definition}[Weakly regular atlas / local \(Q^n\)-action]
\label{def:weakly_regular}
Let \(M^{4n}\) be a compact Hausdorff space. A \emph{weakly regular atlas} (of regularity \(C^r\), \(0\le r\le\infty\)) is an atlas \(\{(U_\alpha^M,\varphi^M_\alpha)\}_{\alpha\in\mathcal A}\) such that for every nonempty overlap \(U_{\alpha\beta}^M\) the overlap map \(\varphi_{\alpha\beta}^M=\varphi_\alpha^M\circ(\varphi_\beta^M)^{-1}\) is a \(C^r\), \(\rho_{\alpha\beta}\)-equivariant diffeomorphism of \(Q^n\)-invariant open subsets of \(\HH^n\) for some automorphism \(\rho_{\alpha\beta}\in Aut(Q^n)\). A \emph{local \(Q^n\)-action} on \(M^{4n}\) is the equivalence class of such atlases under the obvious refinement relation.
\end{definition}

\medskip
\begin{definition}
\label{def:lqtm}
A \emph{local quaternionic toric manifold} is a connected, smooth, $4n$-dimensional manifold $M^{4n}$ equipped with a local $Q^n$-action modeled on the regular representation of $Q^n$ on $\mathbb{H}^n$. 
\end{definition}

\medskip
Given a maximal weakly regular atlas \(\{(U_\alpha^M,\varphi_\alpha^M)\}\) on \(M^{4n}\) define the local orbit projections
\[
\pi_\alpha:\varphi_\alpha^M(U_\alpha^M)\longrightarrow \varphi_\alpha^M(U_\alpha^M)/Q^n.
\]
On overlaps the induced maps are compatible via the \(\rho_{\alpha\beta}\)-equivariant transition maps, hence define an equivalence relation \(\sim_{\mathrm{orb}}\) on the disjoint union \(\coprod_\alpha \varphi_\alpha^M(U_\alpha^M)/Q^n\): For \(b_\alpha\in \varphi_\alpha^M(U_\alpha^M)/Q^n\) and \(b_\beta\in \varphi_\beta^M(U_\beta^M)/Q^n\) we set \(b_\alpha\sim_{\mathrm{orb}} b_\beta\) if \(b_\alpha,b_\beta\) lie in the image of some overlap and are related under the homeomorphism induced by \(\varphi_{\alpha\beta}^M\).

\medskip
\begin{definition}[Orbit space and orbit map]
The \emph{orbit space} \(B_M\) is the quotient
\[
B_M:=\Big(\coprod_\alpha \varphi_\alpha^M(U_\alpha^M)/Q^n\Big)\big/ \sim_{\mathrm{orb}}
\]
equipped with the quotient topology. The natural map
\[
\pi_M:M^{4n}\longrightarrow B_M
\]
induced by the local projections is the \emph{orbit map}.
\end{definition}

The space \(B_M\) is an \(n\)-dimensional manifold with corners, and carries the canonical stratification by codimension:
\[
B_M=\bigsqcup_{k=0}^n \mathcal S^{(k)}B_M,
\]
where \(\mathcal S^{(k)}B_M\) consists of points lying on codimension–\(k\) faces. The strata record isotropy types of the local action: points of \(\mathcal S^{(k)}B_M\) correspond to orbits whose stabilizer is an \((n-k)\)-dimensional quaternionic subtorus.

\subsection{The canonical model}
\label{sec:canonical_model}

In the present subsection, we recall how to construct a $4n$-dimensional manifold equipped with a local $Q^n$-action, from an $n$-dimensional manifold with corners and some linear data of the set of facets.

The transition automorphisms \(\{\rho_{\alpha\beta}\}\) form a \v{C}ech \(1\)-cocycle with values in \( Aut(Q^n)\). This cocycle encodes the \emph{twisting} of the local torus charts and leads to the following construction.

\medskip
\begin{definition}[Principal bundle \(P_M\)]
The cocycle \(\{\rho_{\alpha\beta}\}\) defines a principal \( Aut(Q^n)\)-bundle
\[
\pi_{P_M}:P_M\to B_M,
\]
constructed by gluing \(U_\alpha^B\times Aut(Q^n)\) via
\[
(b_\alpha,h_\alpha)\sim_P (b_\beta,h_\beta)\quad\Longleftrightarrow\quad
b_\alpha=b_\beta\in U_{\alpha\beta}^B,\; h_\alpha=\rho_{\alpha\beta}\circ h_\beta.
\]
\end{definition}

\medskip
Associated to \(P_M\) is the lattice bundle \(\Lambda_M\), with typical fiber \(\Lambda_Q\), obtained via the standard representation \( Aut(Q^n)\cong GL(\Lambda_Q)\) on the lattice.

\begin{remark}[Automorphisms of $Q^n$ and lattice preservation]
Let $Q^n=(S^3)^n \cong \mathrm{SU}(2)^n$. We recall some basic facts about automorphisms of $Q^n$ and lattice preservation following \cite[§3.1]{bg25}. Since $\mathrm{SU}(2)$ has no outer automorphisms, every automorphism of $Q^n$ is obtained by combining:

\begin{itemize}
\item inner automorphisms of each $\mathrm{SU}(2)$–factor 
      (given by conjugation),
\item permutations of the $n$ factors,
\item Weyl group symmetries of each factor.
\end{itemize}

Accordingly, the automorphism group is
\[
\mathrm{Aut}(Q^n) \cong (\mathrm{SO}(3))^n \rtimes S_n,
\]
where $\mathrm{SO}(3)\cong \mathrm{Inn}(\mathrm{SU}(2))$ acts on each factor by conjugation and $S_n$ permutes the factors.

In the present setting, however, the transition functions of a weakly regular atlas preserve the integral lattice 
\[
\Lambda_Q = (4\pi\mathbb{Z})^n \subset \mathfrak{q}^n
\]
and the associated primitive rank-one sublattices. Conjugation acts trivially on $\Lambda_Q$, since lattice elements lie in the center under the exponential map. 

Therefore the structure group reduces to the discrete subgroup
\[
\Gamma := \mathrm{Aut}(Q^n,\Lambda_Q),
\]
consisting of lattice-preserving automorphisms. This subgroup is precisely
\[
\Gamma \cong S_n \ltimes (\mathbb{Z}/2\mathbb{Z})^n,
\]
where $(\mathbb{Z}/2\mathbb{Z})^n$ is the Weyl group acting by sign changes on each lattice coordinate and $S_n$ acts by permutation of coordinates. In particular, $\Gamma$ is a discrete subgroup of 
$\mathrm{Aut}(Q^n)$, and coincides with the lattice automorphism group
\[
\mathrm{GL}(\Lambda_Q).
\]
\end{remark}

\medskip
\begin{definition}[Rank–one sublattice bundle \(\mathcal L_M\)]
Over each codimension–one stratum component of \(B_M\) the local model furnishes a distinguished primitive lattice vector generating a rank–one sublattice of \(\Lambda_Q\). These local rank–one sublattices glue (via the \(\rho_{\alpha\beta}\)-action) to define a rank–one sublattice bundle
\[
\pi_{\mathcal L_M}:\mathcal L_M\to \mathcal S^{(n-1)}B_M.
\]
\end{definition}

\medskip
\begin{definition}[Characteristic pair]\label{def: characteristic_pair}
The pair \((P_M,\mathcal L_M)\) is called the \emph{characteristic pair} associated to the local \(Q^n\)-action on \(M^{4n}\).
\end{definition}

Two remarks are essential:

\begin{remark}
   \begin{enumerate}
   \item[1.] The unimodularity condition: at each codimension–\(k\) point of \(B_M\) the corresponding \(n-k\) rank–one lattice generators must form a direct summand of rank \(n-k\) in \(\Lambda_Q\). This unimodularity guarantees the local model is smoothly equivalent to the standard action on \(\HH^n\).
\item[2.] The principal bundle \(P_M\) is precisely the geometric manifestation of the \emph{twisting} of local torus charts; when the cocycle is trivial, \(P_M\) is trivial as well, and the local action globalizes.
    \end{enumerate}
\end{remark}

From the principal \(Aut(Q^n)\)-bundle \(P_M\) one builds the associated quaternionic torus bundle
\[
\pi_Q:Q_{P_M}:=P_M\times_{ Aut(Q^n)} Q^n \longrightarrow B_M,
\]
whose fibers are (noncanonically) identified with \(Q^n\). Using the lattice data \(\mathcal L_M\) (and the unimodularity condition) one obtains, for each stratum \(\mathcal S^{(k)}B_M\), a well-defined rank–\((n-k)\) subtorus bundle
\[
\pi_{Q_{\mathcal S^{(k)}B}}: Q_{\mathcal S^{(k)}B} \longrightarrow \mathcal S^{(k)}B_M,
\]
encoding the isotropy subtori of the local models.

The canonical model is obtained by quotienting the torus bundle \(Q_{P_M}\) by the equivalence relation that collapses fibers according to the isotropy subtori. In particular, letting \(q,q'\in Q_{P_M}\), then if \(\pi_Q(q)\in\mathcal S^{(k)}B_M\), we say that
\[
q\sim q' \quad\Longleftrightarrow\quad \pi_Q(q)=\pi_Q(q') \quad\text{and}\quad q' q^{-1}\in Q_{\mathcal S^{(k)}B}\big|_{\pi_Q(q)}.
\]

\begin{definition}[Canonical model]\label{def:canonical_model}
The \emph{canonical model} associated with the characteristic pair \((P_M,\mathcal L_M)\) is the quotient space
\[
M_{(P_M,\mathcal L_M)} := Q_{P_M}\big/\!\sim,
\]
equipped with the quotient topology and the natural local \(Q^n\)-action induced from left multiplication on \(Q_{P_M}\). The projection \(M_{(P_M,\mathcal L_M)}\to B_M\) is the orbit map of this action, and the orbit space is \(B_M\).
\end{definition}

The following result compares a local quaternionic toric manifold with its canonical model. Recall that the \emph{Euler class} of a local quaternionic toric manifold is a cohomological invariant 
\[
e(M) \in H^1\!\big(B_M; \mathscr S\big),
\]
where \(\mathscr S\) denotes the sheaf of local sections of the orbit projection \(\pi_M : M \to B_M\). It measures the obstruction to the existence of a global section of \(\pi_M\); in particular, \(e(M)=0\) if and only if \(\pi_M\) admits a global section.

\begin{prop}
\begin{enumerate}
  \item[1.] \(M_{(P_M,\mathcal L_M)}\) is a topological \(4n\)-dimensional manifold carrying a canonical local \(Q^n\)-action with orbit space \(B_M\).
  \item[2.] The canonical model is, locally over each chart \(U_\alpha^B\), equivalent to the standard model \(\pi_{\HH^n}^{-1}(U_\alpha^B)\subset\HH^n\).
  \item[3.] If the \v{C}ech cocycle \(\{\rho_{\alpha\beta}\}\) is cohomologous to the trivial cocycle (equivalently \(P_M\) is trivial) and the Euler obstruction vanishes, then the canonical model is globally \(Q^n\)-equivariantly homeomorphic to the original manifold \(M^{4n}\).
\end{enumerate}
\end{prop}

\section{Homotopy and cohomology of local \(Q^n\)-manifolds}
\label{sec:topology1}

In this section we describe the basic topological structure of manifolds
admitting a local quaternionic torus action. We first describe their fundamental group and then develop the cohomological framework based on the Leray--Serre spectral sequence and its cellular counterpart. These results provide the foundation for the generalized cohomological invariants discussed in the next section.

Let \(M^{4n}\) be a compact smooth manifold equipped with a local \(Q^n\)-action, where \(Q^n=(S^3)^n\) is the quaternionic torus. For simplicity, set $M:=M^{4n}$. We also denote the orbit map by 
\[
\pi_M : M \longrightarrow B_M,
\]
with orbit space \(B_M\) being an \(n\)-dimensional manifold with corners. 

We henceforth assume that \(M\) is connected (so that its orbit space \(B_M\) is also connected) and that the Euler class \(e(M)\) vanishes. In this section we investigate some topological invariants associated to \(M\).

\subsection{On the fundamental groups}
\label{sec:fund_gps}

Fix a base point \(b_0\in \mathrm{Int}(B_M)\) in the interior stratum,
and let \(x_0\in \pi_M^{-1}(b_0)\) be a lift.
Since \(Q^n=(S^3)^n\) is simply connected  one has
\[
\pi_1(Q^n)=0, 
\quad \pi_2(Q^n)=0, 
\quad \pi_3(Q^n)\cong \ZZ^n.
\]
In particular, the fibers of \(\pi_M\) are connected and simply connected.

\medskip

\begin{lemma}\label{lem:sections}
The orbit projection \(\pi_M : M \to B_M\) admits local sections
over contractible neighborhoods of \(b_0\).
\end{lemma}

\begin{proof}
By the definition of a local \(Q^n\)-action, neighborhoods of points in \(B_M\)
are modeled on quotients of open subsets of \(\HH^n\) by linear \(Q^n\)-actions. 
In these models a section is provided by the identity element of \(Q^n\),
and compatibility on overlaps follows from the fact that
transition functions lie in \( Aut(Q^n)\), preserving the group structure.
\end{proof}

\medskip
Choose an open set \(U\) containing $b_0$. Since \(\pi_M|_{U} : \pi_M^{-1}(U)\to U\) is  a locally trivial fiber bundle with typical fiber \(F\cong Q^n\),  the standard long exact sequence of homotopy groups for the fibration applies.  In the local models for a weakly regular atlas one indeed obtains such local fibrations with fiber \(F\cong Q^n\).  When we speak below of the l.e.s. for \(\pi_M\) we mean its l.e.s. in such a local trivializing neighborhood, or globally, whenever \(\pi_M\) is a fibration. 

To write down this standard l.e.s. explicitly, let $\partial$ denote the standard connecting homomorphism in the long exact sequence of fundamental groups, obtained by lifting representatives of spheres in $B_M$ to maps into $M$ and taking their boundaries in the fiber $F\cong Q^n$. Let $i_*$ be the map induced by the inclusion $i : F \hookrightarrow M$ of the typical fiber, and $(\pi_M)_*$ be the map induced by the projection $\pi_M : M \to B_M$. The l.e.s. discussed above is then
\[
\cdots \to \pi_{k+1}(U,b_0)\xrightarrow{\partial} \pi_k(F)\xrightarrow{i_*}\pi_k(M,x_0)
\xrightarrow{(\pi_M)_*}\pi_k(B_M,b_0)\xrightarrow{\partial}\pi_{k-1}(F)\to \cdots
\]
for all \(k\ge 1\).

Since \(\pi_1(Q^n)=0\) and \(\pi_2(Q^n)=0\), the sequence simplifies in low degrees and we obtain the following result.

\medskip

\begin{lemma}\label{lem:exact}
There is a short exact sequence
\[
  1 \;\longrightarrow\; \pi_1(M,x_0) 
  \;\overset{(\pi_M)_*}{\longrightarrow}\; \pi_1(B_M,b_0) 
  \;\longrightarrow\; 1.
\]
\end{lemma}

\begin{proof}
Let \(U \subset B_M\) be a contractible neighborhood of \(b_0\). By Lemma~\ref{lem:sections}, the restriction \(\pi_M^{-1}(U) \to U\) is trivial as a bundle with fiber a neighborhood in \(Q^n\). Thus \(\pi_M^{-1}(U)\) is simply connected in the fiber direction. Since the global fibers are connected, so we have \(\pi_0(Q^n)=0\), and also \(\pi_1(Q^n) = 0\), the long exact sequence of the fibration reduces to
\[
  1 \to \pi_1(M,x_0) \to \pi_1(B_M,b_0) \to 1.
\]
This shows that the homomorphism induced by \(\pi_M\) is injective and surjective, hence an isomorphism.
\end{proof}

\medskip

Summarizing the above we have the following result.

\medskip

\begin{prop}\label{prop:pi1}
For a compact connected manifold \(M^{4n}\) with a local \(Q^n\)-action,
the orbit map induces an isomorphism
\begin{equation}\label{fund1 iso}
\pi_1(M) \;\cong\; \pi_1(B_M).
\end{equation}
\end{prop}

\medskip

\begin{remark}[Comparison with the toric case]
In the case of local $T^n$-actions, Yoshida's \cite{Yoshida} proof of the fundamental group isomorphism $\pi_1(X) \cong \pi_1(B_X)$ requires a careful analysis of the kernel of $\pi_1(X) \to \pi_1(B_X)$, which is essentially the image of $\pi_1(T^n) \cong \mathbb{Z}^n$. This requires the construction of sections and the use of commutative diagrams of split short exact sequences.

In contrast, for local quaternionic toric actions the fiber is $Q^n = (S^3)^n$, which is simply connected. Hence $\pi_1(Q^n) = 0$, and the kernel of $\pi_1(M) \to \pi_1(B_M)$ vanishes automatically. This simplifies the argument considerably: once local sections exist, the long exact sequence of fundamental groups reduces immediately to the isomorphism
\eqref{fund1 iso}.
\end{remark}

\subsection{Cohomology groups for local \(Q^n\)-actions}
\label{sec:cohomology}

In this subsection, we study the cohomology groups for local \(Q^n\)-actions. The Leray-Serre spectral sequence for \(\pi_M\) is the natural tool to compute the cohomology of $M$ from that of the base $B_M$ and of the fiber $Q^n$. So we first describe the general theory and the Leray-Serre  spectral sequence associated with \(\pi_M\). Then, we adapt Yoshida's \cite{Yoshida} cellular construction to the quaternionic setting.

\subsubsection{Leray-Serre spectral sequence}
\label{subsec:leray-serre ss}

For completeness we recall the general (possibly twisted) formulation, before restricting to the untwisted  ($e(M)=0$) case which is relevant for global invariants.

\smallskip
Fix \(\mathbb{Z}\) as the coefficient ring. We describe below the Leray-Serre type spectral sequence associated with \(\pi_M\) (applied locally on trivializing neighborhoods), the arising local coefficient systems, and the concrete simplifications that occur in low degrees due to the fact that the typical fiber is \(Q^n\).

Consider the \(n\)-dimensional manifold with corners structure of \(B_M\) \cite[Proposition 2.7]{bg25}, and let \(\{B^{(p)}_M\}_{p\ge0}\) denote the \(p\)-skeleton of a fixed CW-decomposition of \(B_M\). The successive skeleta filter $B_M$ as follows. Let $\mathcal{S}^{(k)}B_M$ denote the codimension-$k$ stratum of the stratification. The decomposition
\[
B_M = \bigsqcup_{k=0}^{n} \mathcal{S}^{(k)}B_M
\]
defines a natural stratification of $B_M$ by strata. To describe the associated filtration we set
\[
{S}_{(k)}: = \bigcup_{j=0}^{k} \mathcal{S}^{(j)}B_M
\]
so that 
\[
\mathcal{S}_{(0)}B_M \subset \mathcal{S}_{(1)}B_M \subset \cdots \subset \mathcal{S}_{(n)}B_M = B_M,
\]

\iffalse
Note that, for manifolds with corners arising as orbit spaces of locally regular actions, the stratification by codimension coincides (up to homeomorphism) with a CW-structure by faces: each $p$-cell $B^{(p)}$ of $B_M$ corresponds to the relative interior of a codimension-$(n-p)$ face. Thus, the constructions of  $\mathcal{S}^{(k)}B_M$ for strata and $B^{(p)}$ for skeleta, describe equivalent filtrations of $B_M$, and the pullback filtration $\{M^{(p)}\}$ of $M$ may be viewed either combinatorially (via cells) or geometrically (via strata). 

Note that, for manifolds with corners arising as orbit spaces of locally regular actions, the stratification by codimension coincides (up to homeomorphism) with a CW–structure by faces: each $p$–cell $B^{(p)}$ of $B_M$ corresponds to the relative interior of a codimension–$(n-p)$ face. 
Accordingly, the filtration by \emph{skeleta} $\{B^{(p)}\}$ and the filtration by \emph{strata} $\{\mathcal{S}_{(k)}B_M\}$ represent equivalent ways of organizing $B_M$; the former is indexed by the \emph{dimension} of the cells, while the latter is indexed by their \emph{codimension} (dual indexings of the same codimension). We make this precise for the sequence.

\fi

\begin{remark}[Relation between strata and skeleta]
Following Quinn's theory of homotopically stratified spaces \cite{quinn} for an $n$--dimensional manifold with corners $B_M$, the codimension--$k$ stratum \(\mathcal{S}^{(k)}B_M\) coincides (up to homeomorphism) with the open layer between two consecutive skeleta of a compatible CW--structure:
\[
\mathcal{S}^{(k)}B_M \;\simeq\; B^{(n-k)}_M \setminus B^{(n-k-1)}_M.
\]
In other words, each stratum consists of the relative interiors of the $(n-k)$--dimensional faces of $B_M$. Hence, the stratification by codimension $k$ of strata $\{\mathcal{S}_{(k)}B_M\}$ and the CW--filtration by skeleta $\{B^{(p)}_M\}$ of dimension $p$, are dual to each other, related by the index correspondence $p = n - k$.    
\end{remark}

Hence, there is defined a pullback filtration $\{M^{(p)}\}$ of $M$ which may be viewed either combinatorially (via the CW–structure) or geometrically (via the stratification). In detail, for each \(p\) we set
\begin{equation}\label{filtration M}
M^{(p)}:=\pi_M^{-1}\big(B^{(p)}_M\big)\subset M.
\end{equation}
 This defines a filtration on $M$ by subcomplexes (\cite[Chapter 4]{mcc00} and \cite[Chapter 5]{Hat02})
\[
\emptyset=M^{(-1)}\subset M^{(0)}\subset M^{(1)}\subset \cdots \subset M^{(n)}=M.
\]

\noindent Each inclusion of skeleta $(M^{(p-1)} \hookrightarrow M^{(p)})$ defines a pair of spaces $(M^{(p)}, M^{(p-1)})$, and we set
\[
E^{p,q} := H^{p+q}\!\big(M^{(p)}, M^{(p-1)}; \mathbb{Z}\big).
\]
for the relative cohomology between two successive skeleta, and we also define
\[
A^{p,q} := H^{p+q}\!\big(M, M^{(p-1)}; \mathbb{Z}\big),
\]
so that $A^{p,q}$ describes classes of cochains on $M$ that vanish on the $(p-1)$–skeleton, while $E^{p,q}$ records the classes of those cochains supported between successive skeleta.

We also consider the homomorphisms \(i,j,k\) that are defined as 
\begin{align*}
i &: A^{p,q} \longrightarrow A^{p+1,q-1}, &
  j &: A^{p,q} \longrightarrow E^{p,q}, &
  k &: E^{p,q} \longrightarrow A^{p+1,q},
\end{align*}
where
\begin{itemize}
  \item \(i\) is induced from the inclusion of pairs 
  \((M, M^{(p-2)}) \subset (M, M^{(p-1)})\);
  \item \(j\) is induced from the inclusion 
  \((M^{(p)}, M^{(p-1)}) \subset (M, M^{(p-1)})\); and
  \item \(k\) is the connecting homomorphism in the long exact sequence of the triple \((M, M^{(p)}, M^{(p-1)})\).
\end{itemize}
Following the standard construction (see, e.g., \cite[p.37]{mcc00} and \cite{Yoshida} for the toric case), the maps \(i\), \(j\), and \(k\) are precisely the maps appearing in the long exact sequences arranged in the so-called \emph{staircase diagram} \cite[Chapter 5]{Hat02}, in which the horizontal maps correspond to the morphisms $j$ and $k$, while the vertical ones correspond to $i$.
\[
\begin{tikzcd}[column sep=2.8em, row sep=2.8em]
  & \vdots
  \arrow[d,"i"]
  & 
  & \vdots
  \arrow[d,"i"] \\
\cdots 
  \arrow[r,"k"] 
  & A^{p,q} 
      \arrow[r,"j"] 
      \arrow[d,"i"] 
  & E^{p,q} 
      \arrow[r,"k"]  
  & A^{p+1,q} 
      \arrow[r,"j"] 
      \arrow[d,"i"]
  & \cdots \\[-0.2em]
\cdots
  \arrow[r]
  & A^{p+1,q-1} 
      \arrow[r,"j"'] 
      \arrow[d,"i"]
  & E^{p+1,q-1} 
      \arrow[r,"k"'] 
  & A^{p+2,q-1} 
  \arrow[r]
  \arrow[d,"i"]
  & \cdots \\
  & \vdots
  & 
  & \vdots
\end{tikzcd}
\]

\noindent
This staircase diagram describes the long exact sequences obtained from the filtration of $M$ by skeleta. Each diagonal block (consisting a triplet of maps $i,j,k$ making composition meaningful) forms a short exact triangle, and these triangles assemble globally into the diagram above.

\medskip

Before passing to the compressed form of this diagram, we recall the general notion of a spectral sequence (see \cite[Def.~1.1]{mcc00}).

\begin{definition}[Spectral sequence]
A \emph{(first-quadrant, cohomological) spectral sequence} is a sequence of differential bigraded groups
\[
\{(E_r^{p,q}, d_r)\}_{r\ge1}, \qquad p,q\ge0,
\]
where each $E_r^{p,q}$ is a bigraded group and
\[
d_r : E_r^{p,q} \longrightarrow E_r^{p+r,\,q-r+1}
\]
is a differential satisfying $d_r^2 = 0$. The $(r+1)$-page is obtained as the cohomology of $(E_r^{\bullet,\bullet},d_r)$, namely
\[
E_{r+1}^{p,q} \;\cong\;
\ker(d_r : E_r^{p,q}\!\to\!E_r^{p+r,\,q-r+1})
\Big/
\operatorname{im}(d_r : E_r^{p-r,\,q+r-1}\!\to\!E_r^{p,q}).
\]
We call $E_r^{p,q}$ the \emph{$r$-th term} (or \emph{$r$-th page}) of the spectral sequence.
\end{definition}

\noindent Applied to our setting above, we obtain the corresponding {\em exact couple},

\[
\begin{tikzcd}[column sep=large, row sep=large]
\bigoplus_{p,q} A^{p,q} \arrow[rr,"i"] && \bigoplus_{p,q} A^{p,q} \arrow[dl,"j"] \\
& \bigoplus_{p,q} E^{p,q} \arrow[ul,"k"] &
\end{tikzcd}
\]

\noindent
This triangle is exact at each of its three vertices and encodes the same information as the staircase above in a compressed form. An exact couple of this type produces the spectral sequence $\{E^{p,q}_r,d_r\}_{r\geq1}$, where the first page and differential are given by
\[
E_1^{p,q}=E^{p,q}, \quad d_1=j\circ k:E^{p,q}\to E^{p+1,q}.
\]

\noindent
This provides the cohomological spectral sequence associated with the filtration of $M$,  converging to $H^\bullet(M;\mathbb{Z})$.

\medskip

We next identify the local system of coefficients appearing on the $E_2$-page.

\begin{lemma}[Local system associated to a quaternionic fibration]
\label{lem:local_system}
Let \(\pi_M:M\to B_M\) be the orbit projection of a manifold with a local \(Q^n\)-action. Then the collection of groups \(\{H^q(\pi_M^{-1}(b);\mathbb{Z})\}_{b\in B_M}\), together with the induced maps \(\rho_{\alpha\beta}^*\) on overlaps, defines a 
local system of coefficients
\[
\mathscr H^q\;:=\;H^q(Q^n;\mathbb{Z})
\]
on \(B_M\). For every open set \(U\subset B_M\) over which \(\pi_M\) trivializes, \(\mathscr H^q(U)\cong H^q(\pi_M^{-1}(U);\mathbb{Z})\), and on overlaps the identifications are given by the action of \(\mathrm{Aut}(Q^n)\) on \(H^q(Q^n;\mathbb{Z})\).
\end{lemma}

\begin{proof}
Let \(\{(U_\alpha, \varphi_\alpha)\}\) be a weakly regular atlas of local trivializations \(\varphi_\alpha : \pi_M^{-1}(U_\alpha) \xrightarrow{\ \cong\ } U_\alpha \times Q^n\) (see definition \ref{def:weakly_regular}). On overlaps \(U_{\alpha\beta} = U_\alpha \cap U_\beta\),
the transition functions \(\rho_{\alpha\beta} : U_{\alpha\beta} \to \mathrm{Aut}(Q^n)\)
satisfy the \v{C}ech cocycle condition
\[
\rho_{\alpha\beta}\circ\rho_{\beta\gamma}=\rho_{\alpha\gamma.}
\]
By applying the cohomology functor \(H^q(-;\mathbb{Z})\) to these maps, we have automorphisms $\rho^*_{\alpha\beta}: H^q(Q^n;\mathbb{Z})\to H^q(Q^n;\mathbb{Z})$ satisfying the same condition
\[
\rho_{\alpha\beta}^*\circ\rho_{\beta\gamma}^*=\rho_{\alpha\gamma}^*,
\]
so that the cohomology groups of the fibers glue consistently across overlaps. Thus we obtain a cocycle with values in \(\mathrm{Aut}(H^q(Q^n;\mathbb{Z}))\). This cocycle prescribes the gluing of the constant presheaf \(U\mapsto H^q(Q^n;\mathbb{Z})\) across intersections, yielding the desired
locally constant sheaf \(\mathscr H^q\).
\end{proof}

\smallskip
Working on a contractible neighborhood \(U\subset B_M\) where \(\pi_M^{-1}(U)\cong U\times Q^n\), the projection \(\pi_M|_U : \pi_M^{-1}(U)\to U\) is a trivial fibration. The Leray-Serre spectral sequence for this local product has \(E_2\)-term equal to \(H^p(U;H^q(Q^n;\mathbb{Z}))\). The transition automorphisms \(\{\rho_{\alpha\beta}\}\subset \mathrm{Aut}(Q^n)\) of a weakly regular atlas act on \(H^q(Q^n;\mathbb{Z})\), and gluing these local sequences via the Čech machinery (see Lemma \ref{lem:local_system}) produces a global spectral sequence with local coefficients in the associated system \(\mathscr{H}^q\).

\begin{thm}[Leray-Serre spectral sequence for \(\pi_M\)]\label{thm:leray} 
 There is a first-quadrant spectral sequence \((E_r^{p,q},d_r)\) associated with the filtration \eqref{filtration M} of $M$, whose second page is
\[
E_2^{p,q} \;=\; H^p\!\big(B_M;\mathscr H^q\big),
\]
and which converges to the cohomology of $M$ in the sense that
\[
E_\infty^{p,q}
\;\cong\;
\mathrm{gr}^p H^{p+q}(M;\mathbb{Z})
\;:=\; F^p H^{p+q}(M;\mathbb{Z})
  \big/ F^{p+1}H^{p+q}(M;\mathbb{Z}),
\]
where $F^\bullet$ is the filtration on $H^\bullet(M;\mathbb{Z})$ induced by the filtration \eqref{filtration M} of $M$.
\end{thm}

\begin{proof}
We first recall that the orbit projection \(\pi_M : M \to B_M\) is locally modeled on the quotient of an open subset of \(\mathbb{H}^n\) by the left multiplication of the quaternionic torus \(Q^n = (S^3)^n\). By definition \ref{def:weakly_regular} (of a local \(Q^n\)-action), there exists a weakly regular atlas 
\(\{(U_\alpha^M,\varphi_\alpha^M)\}_{\alpha\in\mathcal{A}}\) 
such that each local chart \(\pi_\alpha : 
\pi_M^{-1}(U_\alpha^B) \to U_\alpha^B\) is 
\(Q^n\)-equivariantly diffeomorphic to the trivial bundle
\(U_\alpha^B \times Q^n \to U_\alpha^B\).

\smallskip

On each contractible open subset \(U_\alpha^B \subset B_M\), the classical Leray-Serre spectral sequence for the fibration \(\pi_M^{-1}(U_\alpha^B)\to U_\alpha^B\) applies, since the fiber \(Q^n\) is a compact CW complex.  
This local spectral sequence has $E_2$-term
\[
E_2^{p,q}(U_\alpha^B)
   = H^p\!\big(U_\alpha^B; H^q(Q^n;\mathbb{Z})\big),
\]
and it converges to the cohomology of the total space over $U_\alpha^B$, i.e.
\[
E_\infty^{p,q}(U_\alpha^B) \;\cong\;
\mathrm{gr}^p H^{p+q}\!\big(\pi_M^{-1}(U_\alpha^B);\mathbb{Z}\big).
\]

\smallskip

By Lemma \ref{lem:local_system}, on the overlaps \(U_{\alpha\beta}^B = U_\alpha^B \cap U_\beta^B\), the transition maps of the local trivializations are determined by the automorphisms  \(\rho_{\alpha\beta}\in\mathrm{Aut}(Q^n)\), which induce homomorphisms on fiber cohomology
\[
\rho_{\alpha\beta}^*: H^q(Q^n;\mathbb{Z})
  \longrightarrow H^q(Q^n;\mathbb{Z}).
\]
These maps define an action of the Čech $1$–cocycle \(\{\rho_{\alpha\beta}\}\) on the system of coefficient groups \(H^q(Q^n;\mathbb{Z})\), and hence determine a local coefficient system \(\mathscr{H}^q\) over \(B_M\) with stalk \(\mathscr{H}^q(b) = H^q(\pi_M^{-1}(b);\mathbb{Z})\) for each \(b\in B_M\).

The standard Leray construction \cite[Proposition~15.9]{Switzer} glues the local spectral sequences \(\{E_2^{p,q}(U_\alpha^B)\}\) over the open cover \(\{U_\alpha^B\}\) using this local system of coefficients. The resulting global spectral sequence has $E_2$-term
\[
E_2^{p,q} = H^p(B_M; \mathscr{H}^q)
\]
and converges to the cohomology of the total space \(H^{p+q}(M;\mathbb{Z})\). Finally, since $B_M$ is a compact $n$–manifold with corners, it admits a finite CW decomposition, and the filtration by skeleta
\[
\emptyset = M^{(-1)} \subset M^{(0)} \subset \cdots \subset M^{(n)} = M
\]
is finite. 

The spectral sequence obtained by the exact couple above induces a filtration on the total cohomology $H^r(M;\mathbb{Z})$. 
Specifically, following \cite[p.525]{Hat02} for each~$r$, one defines
\[
F^pH^{r}(M;\mathbb{Z})
:= \operatorname{Im}\!\big(H^{r}(M,M^{(p-1)};\mathbb{Z})
\longrightarrow H^{r}(M;\mathbb{Z})\big),
\]
where on the right is the image of the map induced by the inclusion of pairs $(M,M^{(p-1)})\hookrightarrow (M,\varnothing)$. This filtration measures how cohomology classes of~$M$ ``build up'' as one attaches successive skeleta of~$B_M$ via the orbit projection.

As a result, the spectral sequence produces the natural filtration
\[
0\;=\;F^{r+1}H^{r}(M;\mathbb{Z})
\subset F^{r}H^{r}(M;\mathbb{Z})
\subset\cdots\subset
F^{0}H^{r}(M;\mathbb{Z})=H^{r}(M;\mathbb{Z}),
\]
where $F^pH^{p+q}(M;\mathbb{Z})$ consists of classes whose representatives
are supported on the preimage of the $p$--skeleton of~$B_M$. Hence the spectral sequence lies in the first quadrant and converges strongly to the graded group associated with the filtration of \(H^\bullet(M;\mathbb{Z})\),
\[
\mathrm{gr}^p H^{p+q}(M;\mathbb{Z})
   = F^p H^{p+q}(M;\mathbb{Z})/F^{p+1}H^{p+q}(M;\mathbb{Z}).
\]
This completes the proof.
\end{proof}

We will frequently use the cohomology of the fiber \(Q^n=(S^3)^n\) with integer coefficients. The Künneth formula gives
\[
H^\bullet(Q^n;\mathbb{Z}) \;\cong\; \bigotimes_{i=1}^n H^\bullet(S^3;\mathbb{Z}),
\]
hence the only nonzero cohomology groups in low degree are
\[
H^0(Q^n;\mathbb{Z})\cong\mathbb{Z},\qquad H^3(Q^n;\mathbb{Z})\cong \mathbb{Z}^n,
\]
and higher-degree groups appear in degrees that are sums of multiples of \(3\). Using Theorem~\ref{thm:leray} together with the cohomology of \(Q^n\) one obtains immediate simplifications for the \(E_2\)-page:
\[
E_2^{p,0} \;=\; H^p(B_M;\mathscr H^0) \cong H^p(B_M;\mathbb{Z}),
\qquad
E_2^{p,1}=0,\qquad E_2^{p,2}=0,
\qquad
E_2^{p,3}=H^p\big(B_M;\mathscr H^3\big),
\]
where \(\mathscr H^0\) is the constant local system \(\mathbb{Z}\) and \(\mathscr H^3\) is the rank-\(n\) local system with fibre \(H^3(Q^n;\mathbb{Z})\cong\mathbb{Z}^n\) (the action of \(Aut(Q^n)\) on this lattice gives its monodromy).

\begin{cor}\label{cor:H1H2}
Assume that the local coefficient systems \(\mathscr H^1\) and \(\mathscr H^2\) vanish. Then the following hold:
\begin{enumerate}
  \item[1.] \(\;E_2^{0,1}=E_2^{1,1}=0\) and hence \(H^1(M;\mathbb{Z})\cong H^1(B_M;\mathbb{Z})\).
  \item[2.] Since \(E_2^{p,2}=0\) for all \(p\), the edge-map induces a surjection
    \[
      H^2(M;\mathbb{Z}) \twoheadrightarrow H^2(B_M;\mathbb{Z})
    \]
    In particular, if the Euler class $e(M)$ vanishes, then the edge map is an isomorphism.
\end{enumerate}
\end{cor}

\begin{proof}
From the Leray--Serre spectral sequence \(E_2^{p,q}=H^p(B_M;\mathscr H^q)\) the vanishing of \(\mathscr H^1,\mathscr H^2\) implies there are no nonzero \(E_r^{p,q}\) contributing to total degree \(1\) except \(E_2^{1,0}=H^1(B_M;\mathbb{Z})\). Therefore
\[
H^1(M;\mathbb Z)\;\cong\;H^1(B_M;\mathbb Z).
\]

For total degree \(2\), the only possibly nonzero terms are  \(E_2^{2,0}=H^2(B_M;\mathbb Z)\) and \(E_2^{0,2}=H^0(B_M;\mathscr H^2)=0\), so the edge map yields a surjection
\[
H^2(M;\mathbb Z)\twoheadrightarrow H^2(B_M;\mathbb Z).
\]
If \(\pi_M\) admits a global section so that the bundle is topologically trivial, the sequence collapses at \(E_2\), 
and the edge morphism becomes an isomorphism.
\end{proof}

\medskip
\noindent
The previous corollary established the vanishing of the first and second cohomology groups in algebraic terms via the Leray--Serre spectral sequence. 
We now provide a purely geometric explanation of this phenomenon, following a CW--filtration argument originally due to Scott~\cite[Lemma~3.1]{scott1995}, which applies to local quaternionic toric actions.

\begin{prop}[2--connectedness]
\label{prop:2connected}
Let $M^{4n}$ be a manifold with a local quaternionic torus $Q^n$-action and orbit space $B_M$. If the orbit space $B_M$ is simply connected, then $M$ is $2$--connected; that is,
\[
\pi_1(M)=\pi_2(M)=0.
\]
\end{prop}

\begin{proof}
The argument follows the CW--filtration technique of Scott~\cite[Lemma~3.1]{scott1995}. Each face $F$ of $B_M$ determines a stratum $\pi_M^{-1}(\mathrm{int}\,F)$ modeled on $\mathbb{R}^{4\dim F}\times (S^3)^{n-\dim F}$. By attaching these strata in increasing dimension, we obtain a CW--filtration of $M$ whose $3$--skeleton consists of the $3$--skeleton of $B_M$ together with a wedge of $S^3$--fibers above the vertices of $B_M$. Therefore, by Whitehead’s lemma and results of Milnor (cf.~\cite[Thm.~3.5]{MilnorCW}), the resulting space is $2$--connected.
\end{proof}

\begin{remark}
This 2-connectedness immediately implies $H^1(M;\mathbb{Z})=H^2(M;\mathbb{Z})=0$, in consistency with Corollary \ref{cor:H1H2}.
\end{remark}

\medskip
From the exactness of the spectral sequence in total degree \(p+q=3\) we obtain an extension describing \(\pi_3\)-contributions:  Since \(H^1(M;\mathbb{Z})=H^2(M;\mathbb{Z})=0\), the only nontrivial terms with \(p+q=3\) occur for \((p,q)=(0,3)\) and \((p,q)=(3,0)\) to \(H^3(M;\mathbb{Z})\). Concretely, \(E_2^{p,3}=H^p(B_M;\mathscr H^3)\) and abuts to the graded pieces of \(H^{p+3}(M;\mathbb{Z})\). The first nontrivial piece in total degree \(3\) gives the short exact sequence
\[
0 \longrightarrow E_\infty^{0,3} \longrightarrow H^3(M;\mathbb{Z})
\longrightarrow E_\infty^{3,0} \longrightarrow 0,
\]
with \(E_\infty^{0,3}\) being a quotient of \(E_2^{0,3}=H^0(B_M;\mathscr H^3)\cong \mathbb{Z}^n\) and \(E_\infty^{3,0}\) being a subgroup of \(H^3(B_M;\mathbb{Z})\).

\medskip

So we obtain the following observation.

\medskip
\begin{cor}\label{cor:H3}
If \(B_M\) satisfies \(\pi_4(B_M)=0\) locally (or if the connecting maps \(\partial:\pi_4(B_M)\to\pi_3(Q^n)\) vanish), then one obtains a split exact sequence
\[
0\longrightarrow H^0(B_M;\mathscr H^3)\cong \mathbb{Z}^n
\longrightarrow H^3(M;\mathbb{Z})
\longrightarrow H^3(B_M;\mathbb{Z})\longrightarrow 0,
\]
so \(H^3(M;\mathbb{Z})\) contains a direct summand \(\mathbb{Z}^n\) coming from the fiber.
\end{cor}

%The following appear in Yoshida and here we generalize them

\subsubsection{Cellular spectral sequence}
\medskip
We now adapt Yoshida’s \cite{Yoshida} cellular construction to the quaternionic setting. Let $(M^{4n},Q^n)$ be a manifold with a local $Q^n$-action and orbit map $\pi_M:M\to B_M$. Recall from Section~\ref{sec:canonical_model} that, associated with the characteristic pair 
$(P_M,\mathcal{L}_M)$ on the base manifold with corners $B_M$, there is a canonical 
\emph{characteristic bundle}
\[
\pi_{Q_M} : Q_M \longrightarrow B_M,
\]
a locally trivial quaternionic torus bundle with fiber $Q^n$ and structure group  $\mathrm{Aut}(Q^n)$ (see~\cite[§3.3]{bg25} for details).  

\medskip

Let $B_M$ be equipped with a CW--structure such that each $p$--cell $e^{(p)}$ of $B_M$ is contained in some $k$--dimensional stratum $\mathcal{S}^{(k)}B_M$. We then define the induced filtration on $Q_M$ by
\[
Q_M^{(p)} := \pi_{Q_M}^{-1}\!\big(B_M^{(p)}\big),
\]
where $B_M^{(p)}$ denotes the $p$--skeleton of $B_M$. For each codimension--$k$ stratum $\mathcal{S}^{(k)}B_M \subset B_M$, there exists a rank--$(n-k)$ subbundle
\[
\pi_{S_M} : S_M \longrightarrow \mathcal{S}^{(k)}B_M
\]
of the characteristic bundle $\pi_{Q_M}:Q_M\to B_M$, defined by restricting $\pi_{Q_M}$ to a sufficiently small neighborhood of $\mathcal{S}^{(k)}B_M$ and taking the subbundle whose fibers coincide with the $(n-k)$--dimensional stabilizer tori of the local $Q^n$--action. With this notation, we obtain the following result.

\begin{lemma}\label{restriction}
The restriction of the orbit map $\pi_M:M\to B_M$ to $\mathcal{S}^{(k)}B_M$ is identified with the quotient bundle
\[
Q_M / S_M \;\longrightarrow\; \mathcal{S}^{(k)}B_M.
\]
\end{lemma}

\begin{proof}
The argument follows as in the toric case: let $\tau : D^p \to B_M$ denote the characteristic map of the $p$--cell $e^{(p)} \subset B_M$, where $D^p \subset \mathbb{R}^p$ is the closed unit ball and $\partial D^p$ its boundary sphere. The center $c^{(p)} = \tau(0)$ is the barycenter of $e^{(p)}$, and $\tau(D^p \setminus \partial D^p)$ identifies with the interior of the cell. We denote by
\begin{equation}\label{restriction1}
\nu_{c^{(p)}} : \pi_{Q_M}^{-1}(c^{(p)}) \longrightarrow \pi_M^{-1}(c^{(p)})
\end{equation}
the restriction of the natural projection $\nu : Q_M \to M$ to the fiber over $c^{(p)}$. Then
 there exists a bundle map
\[
\widetilde{\tau}_{Q}:(D^p,\partial D^p)\times \pi_{Q_M}^{-1}(c^{(p)})
\;\longrightarrow\; (Q_M^{(p)},Q_M^{(p-1)}).
\]
That is, the composition $\pi_{Q_M} \circ \widetilde{\tau}_Q$ equals $\tau \circ \mathrm{pr}_{D^p}$ on $D^p \times \pi_{Q_M}^{-1}(c^{(p)})$; in other words, $\widetilde{\tau}_Q$ is a bundle map \emph{covering} $\tau$, and therefore provides a local trivialization of $\pi_{Q_M}$ along the image of $\tau$. Moreover, $\widetilde{\tau}_Q$ restricts to the identity on $\{0\}\times \pi_{Q_M}^{-1}(c^{(p)})$, and it maps $D^p\times S_M(c^{(p)})$ into $S_M$.

This induces a continuous map
\[
\widetilde{\tau}_M:(D^p,\partial D^p)\times \pi_M^{-1}(c^{(p)})
\;\longrightarrow\; (M^{(p)},M^{(p-1)})
\]
which restricts to the identity on $\{0\}\times \pi_M^{-1}(c^{(p)})$. We thus get a diagram

\[
\begin{tikzcd}[column sep=large, row sep=large]
D^p \times \pi_{Q_M}^{-1}(c^{(p)})
  \arrow[r, "id_{D^p} \times \nu_{c^{(p)}}"]
  \arrow[d, "\widetilde{\tau}_Q"']
  \arrow[dr, "\widetilde{\tau}_M" description]
  & D^p \times \pi_M^{-1}(c^{(p)}_\lambda)
      \arrow[d, "\widetilde{\tau}_M"] \\[0.5em]
Q_M
  \arrow[r, "\tau"']
  \arrow[d, "\pi_{Q_M}"']
  & M
      \arrow[d, "\pi_M"] \\[0.5em]
B_M
  \arrow[r, equal]
  & B_M
\end{tikzcd}
\]

 \noindent whose commutativity follows from the fact that $\widetilde{\tau}_M$ is induced from $\widetilde{\tau}_Q$ by the quotient construction $M = Q_M/S_M$. Since both $\pi_{Q_M}$ and $\pi_M$ are projections compatible with this quotient, the naturality of the construction ensures that every square in the diagram  commutes. The construction expresses the compatibility between the local models of $Q_M$ and $M$ over $B_M$.
\end{proof}

We denote by 
\[
C^p_{\!Q}(B_M; \mathscr{H}^q_M)
\subset C^p(B_M; \mathscr{H}^q_M)
\]
the subgroup consisting of those cochains which, on each $p$--cell $e^{(p)} \subset B_M$, take values in the image
\[
\nu_{c^{(p)}}^*\!\big(H^q(\pi_M^{-1}(c^{(p)});\mathbb{Z})\big)
\;\subset\;
H^q(\pi_{Q_M}^{-1}(c^{(p)});\mathbb{Z}),
\]
of the homomorphism 
$\nu_{c^{(p)}}^*$ induced by \eqref{restriction1}. Therefore \(C^p_{\!Q}(B_M; \mathscr{H}^q_M)\) is the subcomplex of the full cochain complex \(C^\bullet(B_M; \mathscr{H}^q_M)\) obtained by restricting to cochains that are compatible with the characteristic bundle $\pi_{Q_M}$.

\medskip

Let $\{(E_Q)^{p,q}_r,d^Q_r\}$ denote the Leray--Serre spectral sequence of the quaternionic torus bundle $\pi_{Q_M}:Q_M\to B_M$. Its construction is completely analogous to that of the spectral sequence $\{(E_M)^{p,q}_r,d^M_r\}$ associated with the orbit projection $\pi_M:M\to B_M$ discussed in  Section~\ref{subsec:leray-serre ss}: using the same filtration of $B_M$ by skeleta, the only difference lies in the replacement of the total space $M$ by the characteristic bundle $Q_M$. In particular, the $E_1$--terms are computed from the cellular cochains $C^p_{\!Q}(B_M; \mathscr{H}^q_M)$ defined above, and the differentials $d^Q_r$ are induced by the attaching maps of the CW--structure on $B_M$.

\medskip

In the next step, we identify the $E_1$-term of this spectral sequence:
Using the bundle maps $\widetilde{\tau}_Q$ introduced above, together with the excision isomorphism and the Künneth formula, we obtain canonical  identifications of the form
\begin{align*}
(E_Q)^{p,q}_1
  &\;\cong\;
  H^{p+q}\!\big((Q_M)^{(p)},(Q_M)^{(p-1)};\mathbb{Z}\big)\\
 & \;\cong\;
  \bigoplus_{\text{$p$-cells } e^{(p)} \subset B_M}
  H^{p}(D^p,\partial D^p;\mathbb{Z})
  \otimes
  H^{q}\!\big(\pi_{Q_M}^{-1}(c^{(p)});\mathbb{Z}\big)    
\end{align*}
where for each $p$--cell $e^{(p)}$ of $B_M$ we fix a characteristic map $\tau:(D^p,\partial D^p)\to(B_M^{(p)},B_M^{(p-1)})$ with barycentric center $c^{(p)}=\tau(0)$. The first isomorphism follows from excision, which replaces each pair $\big((Q_M)^{(p)},(Q_M)^{(p-1)}\big)$ by the disjoint union of the local models $(D^p,\partial D^p)\times \pi_{Q_M}^{-1}(c^{(p)})$ determined by $\tau$. The second isomorphism arises from the Künneth formula applied to these product pairs, using the local triviality of $\pi_{Q_M}$ over the image of 
each $p$--cell of $B_M$.

\iffalse

Similarly, for the orbit projection $\pi_M:M\to B_M$, we have
\begin{align*}
(E_M)^{p,q}_1
  \;\cong\;
  \bigoplus_{\text{$p$-cells } e^{(p)} \subset B_M}
  H^p(D^p,\partial D^p;\mathbb{Z})
  \otimes
  H^q\!\big(\pi_M^{-1}(c^{(p)});\mathbb{Z}\big).
\end{align*}

\begin{remark}
    The two expressions differ only in the total space whose fibers contribute to the coefficient groups $H^q(-;\mathbb{Z})$.
\end{remark}

By the naturality of these identifications and the commutativity of the diagram in the proof of Lemma \ref{restriction}, we get that the following diagram also commutes:
\[
\begin{tikzcd}[column sep=large, row sep=large]
(E_M)^{p,q}_1 
  \arrow[r, "\cong"] 
  \arrow[d] 
  & \displaystyle\bigoplus_\lambda 
    H^p(D^p_\lambda,\partial D^p_\lambda)\otimes H^q(\pi_M^{-1}(c^{(p)}_\lambda))
    \arrow[d] \\
(E_Q)^{p,q}_1 
  \arrow[r, "\cong"] 
  & \displaystyle\bigoplus_\lambda 
    H^p(D^p_\lambda,\partial D^p_\lambda)\otimes H^q(\pi_{Q_M}^{-1}(c^{(p)}_\lambda))
\end{tikzcd}
\]

Moreover, the induced homomorphism
\[
\nu^*_{c^{(p)}_\lambda}:H^\bullet(\pi_M^{-1}(c^{(p)}_\lambda);\mathbb{Z})\to
H^\bullet(\pi_{Q_M}^{-1}(c^{(p)}_\lambda);\mathbb{Z})
\]
is injective. Therefore
\[
(E_M)^{p,q}_1 \;\cong\; C^p(B_M;\mathscr H^q_M),
\]
where $\mathscr H^q_M$ is the local system with fiber $H^q(Q^n;\mathbb{Z})$ determined by the characteristic bundle. Thus the $E_1$-page is identified with the cellular cochains of $B_M$ with coefficients in $\mathscr H^q_M$.

\medskip
Summarizing all the above we have the following result.

\fi

\medskip
\begin{thm}\label{thm:cellularSSS}
For a local quaternionic toric action $\pi_M:M\to B_M$, the cohomology spectral sequence $\{(E_M)^{p,q}_r,d_r\}$ satisfies
\[
(E_M)^{p,q}_1 \;\cong\; C^p(B_M;\mathscr H^q_M), \qquad
(E_M)^{p,q}_2 \;\cong\; H^p(B_M;\mathscr H^q_M).
\]
\end{thm}

\begin{proof}
We construct the spectral sequence associated with the filtration of $M$ by the preimages of the skeleta of $B_M$:
\[
M^{(p)} := \pi_M^{-1}(B_M^{(p)}), \qquad p\ge0.
\]
The corresponding cohomological spectral sequence is
\[
(E_M)^{p,q}_r \;=\; H^{p+q}(M^{(p)},M^{(p-1)};\mathbb{Z}),
\]
with differentials $d_r : (E_M)^{p,q}_r \to (E_M)^{p+r,q-r+1}_r$. We first describe the $E_1$--term.  
For each $p$--cell $e^{(p)} \subset B_M$, let 
\[
\tau : (D^p,\partial D^p) \longrightarrow (B_M^{(p)},B_M^{(p-1)})
\]
be a characteristic map with barycentric center $c^{(p)} = \tau(0)$.  
By excision, the relative cohomology of the pair $\big(M^{(p)},M^{(p-1)}\big)$ is identified with the direct sum over all $p$--cells of the relative cohomologies of the corresponding local models:
\[
H^{p+q}\big(M^{(p)},M^{(p-1)};\mathbb{Z}\big)
\;\cong\;
\bigoplus_{e^{(p)}}
H^{p+q}\big(\pi_M^{-1}(\tau(D^p)),\pi_M^{-1}(\tau(\partial D^p));\mathbb{Z}\big).
\]

Over each $p$--cell, the bundle $\pi_M$ is locally trivial with fiber $\pi_M^{-1}(c^{(p)})$.  
Hence, by the Künneth formula and naturality of the product decomposition, we have canonical isomorphisms
\[
H^{p+q}\big(\pi_M^{-1}(\tau(D^p)),\pi_M^{-1}(\tau(\partial D^p));\mathbb{Z}\big)
\;\cong\;
H^p(D^p,\partial D^p;\mathbb{Z})
\otimes
H^q(\pi_M^{-1}(c^{(p)});\mathbb{Z}).
\]
Summing over all $p$--cells, we obtain
\[
(E_M)^{p,q}_1
\;\cong\;
\bigoplus_{\text{$p$-cells } e^{(p)} \subset B_M}
H^p(D^p,\partial D^p;\mathbb{Z})
\otimes
H^q(\pi_M^{-1}(c^{(p)});\mathbb{Z}).
\]

\smallskip
Similarly, for the characteristic bundle $\pi_{Q_M} : Q_M \to B_M$, we have
\begin{align*}
(E_Q)^{p,q}_1
  \;\cong\;
  \bigoplus_{\text{$p$-cells } e^{(p)} \subset B_M}
  H^p(D^p,\partial D^p;\mathbb{Z})
  \otimes
  H^q\!\big(\pi_{Q_M}^{-1}(c^{(p)});\mathbb{Z}\big).
\end{align*}

Notice that the two expressions for $(E_M)^{p,q}_1$ and $(E_Q)^{p,q}_1$ differ only in the total space whose fibers contribute to the coefficient groups $H^q(-;\mathbb{Z})$.

By the naturality of these identifications and the commutativity of the diagram in Lemma~\ref{restriction}, we obtain the following commutative diagram:
\[
\begin{tikzcd}[column sep=large, row sep=large]
(E_M)^{p,q}_1 
  \arrow[r, "\cong"] 
  \arrow[d] 
  & \displaystyle\bigoplus_{e^{(p)}} 
    H^p(D^p,\partial D^p;\mathbb{Z}) \otimes H^q(\pi_M^{-1}(c^{(p)});\mathbb{Z})
    \arrow[d] \\
(E_Q)^{p,q}_1 
  \arrow[r, "\cong"] 
  & \displaystyle\bigoplus_{e^{(p)}} 
    H^p(D^p,\partial D^p;\mathbb{Z}) \otimes H^q(\pi_{Q_M}^{-1}(c^{(p)});\mathbb{Z})
\end{tikzcd}
\]

Moreover, the induced homomorphisms
\[
\nu^*_{c^{(p)}} : 
H^\bullet(\pi_M^{-1}(c^{(p)});\mathbb{Z})
\longrightarrow
H^\bullet(\pi_{Q_M}^{-1}(c^{(p)});\mathbb{Z})
\]
are injective. Hence, the cochains on $B_M$ taking values in these fibers form a well-defined subcomplex of $C^p(B_M;\mathscr H^q_M)$, where $\mathscr H^q_M$ denotes the local system with stalk $H^q(Q^n;\mathbb{Z})$ determined by the characteristic bundle.

It follows that the $E_1$--page of the spectral sequence can be identified as
\[
(E_M)^{p,q}_1 \;\cong\; C^p(B_M;\mathscr H^q_M),
\]
the cellular cochains of $B_M$ with local coefficients in $\mathscr H^q_M$.

\smallskip

Finally, the differential $d_1$ on the $E_1$--page corresponds precisely to the cellular coboundary operator of this cochain complex.  
Therefore, taking cohomology with respect to $d_1$ yields
\[
(E_M)^{p,q}_2
\;\cong\;
H^p(B_M;\mathscr H^q_M),
\]
the cohomology of $B_M$ with coefficients in the local system $\mathscr H^q_M$.  
This completes the proof.
\end{proof}

\medskip
\begin{remark}
\begin{enumerate}
\item[1.] For $q=0$, $(E_M)^{p,0}_2\cong H^p(B_M;\mathbb{Z})$.
\item[2.] Since $H^1(Q^n;\mathbb{Z})=H^2(Q^n;\mathbb{Z})=0$, it follows that $(E_M)^{p,1}_r=(E_M)^{p,2}_r=0$ for all $p,r$, because the $E_2$--page of the Leray--Serre spectral sequence
is given by
\[
(E_M)^{p,q}_2 \;\cong\; H^p\!\big(B_M;\mathscr{H}^q_M\big),
\]
where $\mathscr{H}^q_M$ is the local system with stalk $H^q(Q^n;\mathbb{Z})$. Hence if $H^q(Q^n;\mathbb{Z})=0$, the entire row $q$ of the spectral sequence vanishes identically at $E_2$, and therefore at all subsequent pages $(E_M)^{p,q}_r$ for $r\ge2$.
\item[3.] If $B_M$ is an oriented surface with $\partial B_M\neq\emptyset$, the spectral sequence degenerates at the $E_2$-term, for the same reason as in the toric case: fibers over $0$-cells in $\partial B_M$ are of dimension $\leq 1$, hence cohomologically trivial in positive degrees.
\end{enumerate}
\end{remark}

\begin{remark}[On exotic terms]
A closely related approach has been recently developed by Royo Prieto and Saralegi-Aranguren \cite{royo-saralegi}, who constructed a \emph{Gysin braid} describing the relationship between the cohomology of a manifold $M$ with a smooth $S^3$-action and the (intersection) cohomology of its orbit space $M/S^3$. Their construction extends the classical Gysin sequence to the nonfree and semifree cases, where $M/S^3$ becomes a stratified pseudomanifold, and additional ``exotic'' terms arise from the fixed-point set of the $S^1 \subset S^3$ subgroup.
Unlike the Gysin braid construction of \cite{royo-saralegi}, no genuinely ``exotic'' terms appear in the spectral sequences considered here. In our setting the orbit space $B_M$ of a locally standard $Q^n$-action is a smooth manifold with corners rather than a pseudomanifold, so intersection cohomology does not intervene. Nevertheless, when the local system $\mathscr{H}^3$ determined by the monodromy of the $Q^n$-action is nontrivial, the induced extensions on the $E_\infty$-page can be viewed as algebraic analogues of such exotic contributions, reflecting the twisting of the cohomology of the fibers along $B_M$.
\end{remark}

\paragraph*{Quoric manifolds.}\label{par:quoric}
In \cite{ho}, Hopkinson introduced \emph{quoric manifolds} as quaternionic analogues of quasitoric manifolds \cite{dj}. Recall that a quoric manifold is a smooth $4n$-dimensional manifold \( M:=M^{4n} \) equipped with a locally regular quaternionic torus $Q^n$-action whose orbit space is an $n$-dimensional simple polytope $P:=P^n$ with projection map $\pi:M\to P$. Hopkinson computes the cohomology of global quoric manifolds by using the classical \emph{Borel construction} (see \cite[\S 6]{ho}), which involves the fibration:
\[
Q^n \longrightarrow EQ^n \times M \longrightarrow B_QM,
\]
where \( B_QM = (EQ^n \times M)/Q^n \) is the total space of the Borel construction. This construction models the cohomology of the manifold \( M \) as the cohomology of the base space \( B_QM \), which is shown to be homotopy equivalent to the quaternionic \emph{Davis-Januszkiewicz space} \( DJ_Q(K_P) \). The cohomology of \( DJ_Q(K_P) \) is computed as a \emph{Stanley-Reisner ring}:
\[
H^\bullet(DJ_Q(K_P)) \cong \mathbb{Z}[u_1, \dots, u_m] / I_P,
\]
where \( u_a \) corresponds to a generator in degree 4, and the relations are given by the ideal \( I_P \) associated with the polytope \( P \). This result establishes a key connection between the cohomology of \( M \) and the combinatorial structure of the underlying polytope.

\noindent

In the following example we see how the cohomology of a quoric manifold can be computed using the spectral sequences developed above (Theorem \ref{thm:cellularSSS}). We now show how the cellular spectral sequence of Theorem~\ref{thm:cellularSSS} recovers the expected answer to the simplest quoric example \(M=\mathbb H P^1\cong S^4\).

\begin{example}[The basic quoric manifold \(\mathbb H P^1\)]
Let the quoric manifold \(M=\mathbb H P^1\cong S^4\) with the standard \(Q\)-action. The orbit space is the one-dimensional simple polytope
\[
P=\Delta^1=[0,1]
\]
with orbit projection
\[
\pi:S^4\longrightarrow [0,1].
\]
Over the interior of the interval the orbit is \(S^3\), while over the two endpoints the orbit collapses to a point. Thus \(\pi\) is a stratified orbit projection.

Following Theorem \ref{thm:cellularSSS}, we compute the cellular spectral sequence associated with the filtration
\[
M^{(p)}=\pi^{-1}(P^{(p)}).
\]
The \(0\)-skeleton of \(P\) consists of two vertices \(v_0,v_1\), and
\[
M^{(0)}=\pi^{-1}(\{v_0,v_1\})
\]
is the two-point fixed set. Hence
\[
E_1^{0,0}\cong H^0(M^{(0)};\mathbb Z)\cong \mathbb Z^2,
\]
and \(E_1^{0,q}=0\) for \(q>0\).

The unique open \(1\)-cell has orbit \(S^3\). Equivalently,
\[
M^{(1)}/M^{(0)}\cong \Sigma S^3\cong S^4.
\]
Thus the only nonzero terms in column \(p=1\) are
\[
E_1^{1,0}\cong \mathbb Z,
\qquad
E_1^{1,3}\cong \mathbb Z.
\]
The differential
\[
d_1:E_1^{0,0}\cong \mathbb Z^2
\longrightarrow
E_1^{1,0}\cong \mathbb Z
\]
is the ordinary cellular coboundary of the interval. With the standard orientation, it is given by
\[
d_1(a,b)=b-a.
\]
Therefore
\[
E_2^{0,0}\cong \mathbb Z,
\qquad
E_2^{1,0}=0.
\]
Since \(E_1^{0,3}=0\), the class in \(E_1^{1,3}\cong \mathbb Z\) survives to \(E_2\). There are no higher differentials because \(P\) is one-dimensional. Hence
\[
E_\infty^{0,0}\cong \mathbb Z,
\qquad
E_\infty^{1,3}\cong \mathbb Z.
\]
Reading off total degree \(k=p+q\), we obtain
\[
H^k(\mathbb H P^1;\mathbb Z)
\cong
\begin{cases}
\mathbb Z, & k=0,4,\\
0, & \text{otherwise}.
\end{cases}
\]
This agrees with the standard computation
\[
H^\ast(S^4;\mathbb Z)\cong \mathbb Z[u]/(u^2),
\qquad |u|=4.
\]
\end{example}

We also recover the case \(M=\mathbb H P^2\).

\begin{example}[The case \(M=\mathbb H P^2\)]
Let \(M=\mathbb H P^2\), equipped with the standard quoric action. The orbit space is the simplex
\[
P=\Delta^2.
\]
Thus \(P\) has three vertices, three edges, and one two-dimensional face. The cellular spectral sequence of Theorem~\ref{thm:cellularSSS} has
\[
E_1^{p,q}
\cong
\bigoplus_{\dim F=p} H^q(Q^F;\mathbb Z),
\qquad
Q^F\cong (S^3)^{\dim F}.
\]
Hence the nonzero \(E_1\)-terms are
\[
E_1^{0,0}\cong \mathbb Z^3,
\qquad
E_1^{1,0}\cong \mathbb Z^3,
\qquad
E_1^{2,0}\cong \mathbb Z,
\]
coming from \(H^0\) of the orbits,
\[
E_1^{1,3}\cong \mathbb Z^3,
\qquad
E_1^{2,3}\cong H^3((S^3)^2;\mathbb Z)\cong \mathbb Z^2,
\]
and
\[
E_1^{2,6}\cong H^6((S^3)^2;\mathbb Z)\cong \mathbb Z.
\]

On the \(q=0\) row, the differential is the ordinary cellular coboundary of the simplex \(\Delta^2\). Since \(\Delta^2\) is contractible, this gives
\[
E_2^{0,0}\cong \mathbb Z,
\qquad
E_2^{1,0}=0,
\qquad
E_2^{2,0}=0.
\]

On the \(q=3\) row, the differential
\[
d_1:E_1^{1,3}\cong \mathbb Z^3
\longrightarrow
E_1^{2,3}\cong \mathbb Z^2
\]
is determined by the three edge inclusions into the unique two-dimensional face, together with the corresponding maps on
\[
H^3(S^3;\mathbb Z)\longrightarrow H^3((S^3)^2;\mathbb Z).
\]
For the standard characteristic data of \(\mathbb H P^2\), this map has rank \(2\).
Consequently
\[
E_2^{1,3}\cong \mathbb Z,
\qquad
E_2^{2,3}=0.
\]
Finally, the class
\[
E_1^{2,6}\cong \mathbb Z
\]
has no possible incoming \(d_1\)-differential and survives to \(E_2\). Since the base is two-dimensional, there are no higher differentials. Therefore
\[
E_\infty^{0,0}\cong \mathbb Z,
\qquad
E_\infty^{1,3}\cong \mathbb Z,
\qquad
E_\infty^{2,6}\cong \mathbb Z.
\]
Reading off total degree \(k=p+q\), we obtain
\[
H^k(\mathbb H P^2;\mathbb Z)
\cong
\begin{cases}
\mathbb Z, & k=0,4,8,\\
0, & \text{otherwise}.
\end{cases}
\]
This agrees with the standard formula
\[
H^\ast(\mathbb H P^2;\mathbb Z)
\cong
\mathbb Z[u]/(u^3),
\qquad |u|=4.
\]
\end{example}

\begin{remark}[The general projective case]
For \(M=\mathbb H P^n\), the orbit space is the simplex \(\Delta^n\). The cellular spectral sequence has
\[
E_1^{p,q}
\cong
\bigoplus_{\dim F=p}
H^q((S^3)^p;\mathbb Z).
\]
Equivalently, since
\[
H^\ast((S^3)^p;\mathbb Z)
\cong
\Lambda_{\mathbb Z}(x_1,\dots,x_p),
\qquad |x_i|=3,
\]
(see \cite[Proposition~1.4]{fung}) the \(E_1\)-page is determined by the face lattice of \(\Delta^n\), with coefficients in exterior algebras attached to the corresponding orbits. 

The differential \(d_1\) is the cellular incidence differential of the face complex, decorated by the coefficient maps coming from the quotient maps between adjacent orbits. If \(F'\subset\partial F\) is a codimension-one face inclusion, then the orbit over \(F'\) is obtained from the orbit over \(F\) by collapsing the \(S^3\)-factor determined by the corresponding characteristic vector. Hence, after choosing coordinates compatible with the characteristic data, the induced map on orbit cohomology is the corresponding homomorphism of exterior algebras, omitting the generator associated with the collapsed factor.

The computations for \(\mathbb H P^1\) and \(\mathbb H P^2\) suggest the general pattern: the spectral sequence recovers the additive cohomology groups
\[
H^k(\mathbb H P^n;\mathbb Z)
\cong
\begin{cases}
\mathbb Z, & k=0,4,\dots,4n,\\
0, & \text{otherwise}.
\end{cases}
\]
Together with Hopkinson's ring computation (see \cite[\S 6.3]{ho}), this is compatible with
\[
H^\ast(\mathbb H P^n;\mathbb Z)
\cong
\mathbb Z[u]/(u^{n+1}),
\qquad |u|=4.
\]
\end{remark}

\subsubsection{Euler characteristic}

We now compute the Euler characteristic of $M$ directly from the first page of the spectral sequence. In the toric case, $S^{(0)}{B_X}$ corresponds to the fixed point set of the $T^n$-action. In the quaternionic case, $S^{(0)}{B_M}$ corresponds to fibers of $\pi_M$ which reduce to a point, i.e.\ $Q^n$-fixed points of the local quaternionic action. Thus the following result verifies the quaternionic analogue of the fact that the Euler characteristic equals the number of fixed points.

\medskip

\begin{cor}[Euler characteristic]
Assume $B_M$ is a finite CW complex. Then the Euler characteristic $\chi(M)$ is equal to the cardinality of the set $S^{(0)}{B_M}$ of $0$-cells of $B_M$.
\end{cor}
\begin{proof}
As in the toric case, using rational coefficients and the identification $(E_M)^{p,q}_1\cong C^p(B_M;\mathscr H^q_M)$ we compute
\[
\chi((E_M)_1)=\sum_{p,q}(-1)^{p+q}\dim(E_M)^{p,q}_1
=\sum_p\sum_\lambda (-1)^p\chi\big(\pi_M^{-1}(c^{(p)}_\lambda)\big).
\]
Each fiber $\pi_M^{-1}(c^{(p)}_\lambda)$ is a quaternionic torus of dimension $\leq n$. Its Euler characteristic vanishes unless the fiber is a point, which happens precisely at $0$-cells. Thus $\chi((E_M)_1)$ is the cardinality of $S^{(0)}{B_M}$. On the other hand, $\chi((E_M)_r)$ is independent of $r$, and $\chi((E_M)_\infty)=\chi(M)$. Hence \(\chi(M)=|S^{(0)}{B_M}|\).
\end{proof}

\begin{remark}
In particular, $\chi(M)$ is determined combinatorially by the orbit space $B_M$.
\end{remark}

\section{$K$- and $L$-Theory Invariants}
\label{sec:topology2}

In this section we study secondary topological invariants of manifolds with a local quaternionic torus action. These include the quaternionic analogues of complex $K$-theory classes, as well as $L$-theory and signature-type invariants arising from the intersection form on middle-dimensional homology. We also discuss Novikov additivity and related structural properties.

\subsection{K-theory invariants}
\label{sec:Ktheory}

In this subsection we explain that the cellular spectral sequence constructed above has a direct analogue in complex \(K\)-theory. This is the sense in which the cohomological construction may be repeated with \(H^\bullet(-)\) replaced by \(K^\bullet(-)\). Since the orbit projection of a local quaternionic toric manifold is stratified, the appropriate construction is the spectral sequence associated with the filtration induced from the skeleta of the orbit space. This is the \(K\)-theoretic analogue of the cellular spectral sequence, in the
same general spirit as the Atiyah--Hirzebruch spectral sequence
\cite{AH61}.

We now give the \(K\)-theoretic analogue of Theorem~\ref{thm:cellularSSS}. We still assume that \(B_M\) is equipped with a CW-complex structure compatible with the orbit-type stratification, so that each \(p\)-cell \(e^{(p)}\) is contained in a single stratum of \(B_M\).

Let
\[
\pi_{Q_M}:Q_M\longrightarrow B_M
\]
be the orbit projection. For each \(q\), let \(\mathscr K^q_{Q_M}\) denote the local system on \(B_M\) whose stalk at \(b\in B_M\) is
\[
\mathscr K^q_{Q_M}(b)
=
K^q\bigl(\pi_{Q_M}^{-1}(b)\bigr).
\]
We denote by
\[
\bigl(C^p(B_M;\mathscr K^q_{Q_M}),\delta\bigr)
\]
the cellular cochain complex of \(B_M\) with coefficients in this local system.

Let \(c^{(p)}\) be the barycenter of a \(p\)-cell \(e^{(p)}\subset B_M\). The natural quotient map
\[
\nu:Q_M\longrightarrow M
\]
restricts over \(c^{(p)}\) to a map
\[
\nu_{c^{(p)}}:\pi_{Q_M}^{-1}(c^{(p)})
\longrightarrow
\pi_M^{-1}(c^{(p)}).
\]
Since complex \(K\)-theory is contravariant, this induces a homomorphism
\[
\nu_{c^{(p)}}^\ast:
K^q\bigl(\pi_M^{-1}(c^{(p)})\bigr)
\longrightarrow
K^q\bigl(\pi_{Q_M}^{-1}(c^{(p)})\bigr).
\]
We define
\[
C^p(B_M;\mathscr K^q_M)
\subset
C^p(B_M;\mathscr K^q_{Q_M})
\]
to be the subgroup consisting of those cellular cochains whose value on each \(p\)-cell \(e^{(p)}\) lies in the image of \(\nu_{c^{(p)}}^\ast\). Equivalently,
\[
C^p(B_M;\mathscr K^q_M)
=
\left\{
\alpha\in C^p(B_M;\mathscr K^q_{Q_M})
\;\middle|\;
\alpha(e^{(p)})\in
\operatorname{Im}\left(
\nu_{c^{(p)}}^\ast:
K^q(\pi_M^{-1}(c^{(p)}))
\to
K^q(\pi_{Q_M}^{-1}(c^{(p)}))
\right)
\right\}
\]
for every $p$-cell $e^{(p)}$.
The differential on \(C^\bullet(B_M;\mathscr K^q_M)\) is induced by the cellular coboundary of \(B_M\), together with the maps in \(K\)-theory determined by the quotient maps between adjacent orbit fibers.

\begin{thm}[Cellular \(K\)-theory spectral sequence] \label{thm:K-theory-SS} 
For a local quaternionic toric action
\[
\pi_M:M\longrightarrow B_M,
\]
the \(K\)-theory spectral sequence associated with the filtration
\[
M^{(p)}=\pi_M^{-1}(B_M^{(p)})
\]
satisfies
\[
(E_M)^{p,q}_1
\cong
C^p(B_M;\mathscr K^q_M),
\qquad
(E_M)^{p,q}_2
\cong
H^p(B_M;\mathscr K^q_M).
\]
The spectral sequence converges to \(K^{p+q}(M)\), with respect to the filtration of \(K^\bullet(M)\) induced by the skeleta of \(B_M\).
\end{thm}

\begin{proof}
The proof is the standard exact-couple construction for a finite filtered space (see \cite[Chapter~2]{mcc00}); the same construction applies to generalized cohomology theories, in particular to complex \(K\)-theory. We apply it to the filtration
\[
M^{(p)}=\pi_M^{-1}(B_M^{(p)}).
\]
Applying complex \(K\)-theory to the relative pairs
\[
(M^{(p)},M^{(p-1)})
\]
gives an exact couple and therefore a spectral sequence with
\[
(E_M)^{p,q}_1
=
K^{p+q}\bigl(M^{(p)},M^{(p-1)}\bigr)
\]
converging to \(K^{p+q}(M)\), with respect to the filtration induced on \(K^\bullet(M)\).

We now identify the first page. Let \(e^{(p)}_\lambda\) be a \(p\)-cell of
\(B_M\), and let
\[
\tau_\lambda:(D^p,\partial D^p)
\longrightarrow
(B_M^{(p)},B_M^{(p-1)})
\]
be its characteristic map. Since the CW-structure is compatible with the orbit-type stratification, the orbit type is constant over the interior of \(e^{(p)}_\lambda\). If \(b^{(p)}_\lambda=\tau_\lambda(0)\) is the barycenter, then, after excision, the corresponding relative piece is identified with
\[
(D^p,\partial D^p)\times \pi_M^{-1}(b^{(p)}_\lambda).
\]
Thus, by excision and the suspension isomorphism in complex \(K\)-theory,
\[
K^{p+q}\bigl(
\pi_M^{-1}(\tau_\lambda(D^p)),
\pi_M^{-1}(\tau_\lambda(\partial D^p))
\bigr)
\cong
K^{p+q}\bigl(
(D^p,\partial D^p)\times \pi_M^{-1}(b^{(p)}_\lambda)
\bigr)
\cong
K^q\bigl(\pi_M^{-1}(b^{(p)}_\lambda)\bigr).
\]
Summing over all \(p\)-cells gives
\[
(E_M)^{p,q}_1
\cong
\bigoplus_\lambda
K^q\bigl(\pi_M^{-1}(b^{(p)}_\lambda)\bigr).
\]
By the definition of the cellular coefficient system \(\mathscr K^q_M\), the right-hand side is precisely
\[
C^p(B_M;\mathscr K^q_M).
\]
Hence
\[
(E_M)^{p,q}_1
\cong
C^p(B_M;\mathscr K^q_M).
\]

Under this identification, the first differential \(d_1\) is induced by the boundary map in the exact couple. Equivalently, it is the cellular coboundary determined by the attaching maps of \(B_M\), together with the maps in \(K\)-theory coming from the quotient maps between adjacent orbit fibers. Therefore, taking cohomology with respect to \(d_1\), we obtain
\[
(E_M)^{p,q}_2
\cong
H^p(B_M;\mathscr K^q_M).
\]
This completes the proof.
\end{proof}

\begin{remark}[Existence and convergence]
The spectral sequence above is an instance of the standard spectral sequence associated with a finite filtered space and a generalized cohomology theory; see \cite[Thm.~4.7]{LueckDavis98}. In the present setting the filtration is
\[
M^{(p)}=\pi_M^{-1}(B_M^{(p)}),
\]
where \(B_M^{(p)}\) denotes the \(p\)-skeleton of a CW-structure on \(B_M\) compatible with the orbit-type stratification. Since \(B_M\) is compact and the filtration is finite, the spectral sequence converges strongly to
\(K^\bullet(M)\).
\end{remark}

\begin{remark}[The untwisted case]
The vanishing of the Euler class \(e(M)\) implies that the orbit projection admits a global section in the sense of the canonical model. Consequently, the twisting coming from the obstruction to the existence of a section does not appear in the cellular spectral sequence. In this sense, the spectral sequence above is the untwisted version of the \(K\)-theoretic cellular spectral sequence.

However, the orbit projection remains stratified: the orbit over a face \(F\subset B_M\) depends on the dimension and isotropy type of \(F\). Therefore the coefficient groups on the \(E_1\)-page are still face-dependent:
\[
(E_M)^{p,q}_1
\cong
\bigoplus_{\dim F=p} K^q(Q^F),
\]
where \(Q^F\) denotes the orbit over the relative interior of \(F\).
\end{remark}

\begin{remark}[Twisted refinements]
If the Euler obstruction does not vanish, one expects a twisted version of the cellular spectral sequence, in which the coefficient systems over the face poset are modified by the obstruction data. This is analogous in spirit to twisted Atiyah--Hirzebruch spectral sequences in complex \(K\)-theory; see \cite{Rosenberg82,AtiyahSegal05}. We do not develop this twisted theory here.
\end{remark}

\subsection{On signatures in the oriented quaternionic four-dimensional case}
\label{sec:sign}

The four-dimensional case is the first nontrivial instance in which intersection theory and signature phenomena appear for local quaternionic toric manifolds. Indeed, the signature of a compact oriented manifold is defined only in real dimensions $4k$, as it arises from the symmetric bilinear intersection form on $H^{2k}(M;\mathbb{Z})$. For $M^{4}$, this reduces to the intersection form on $H^2(M)$, whose signature provides a fundamental topological invariant. In the toric setting, this invariant is computed via the Meyer cocycle associated to the monodromy of the torus bundle. The goal of this section is to construct and compute its quaternionic analogue for four-dimensional manifolds endowed with local quaternionic toric actions. We show that the signature can be obtained through Novikov additivity, with local contributions governed by a quaternionic Meyer cocycle valued in $\mathbb{Z}$.

In this subsection, we describe a method of computing the signature of four-dimensional quaternionic toric manifolds under local quaternionic toric actions, using Novikov additivity.  We assume that both $M$ and the interior of $B_M$ are oriented so that a weakly standard atlas of $M$ used in this subsection and the induced atlas of $B_M$ are compatible with the given orientations.   Here, the assumption $e(M)=0$ is not necessary.

We now compute the global signature $\sigma(M)$ using a decomposition of $B_M$ into interior and boundary parts, applying Novikov additivity to the corresponding preimages in $M$.

For simplicity, suppose that $B_M$ has only one boundary component with $S^{(0)}{B_M}\neq \varnothing$. We divide $B_M$ into two parts $(B_M)_1$ and $(B_M)_2$, where $(B_M)_2$ is the closed neighborhood of the boundary $\partial B_M$ such that $\partial B_M$ is a deformation retract of $(B_M)_2$, and $(B_M)_1$ is the closure $(B_M)_1=(B_M\setminus (B_M)_2)^{\mathrm{cl}}$ of the remainder. We divide $B_M$ into two parts $(B_M)_1$ and $(B_M)_2$, where $(B_M)_2$ is a closed collar neighborhood of the boundary $\partial B_M$ such that $\partial B_M$ is a deformation retract of $(B_M)_2$, and $(B_M)_1$ is the closure of the complement, 
\[
(B_M)_1 = \overline{B_M \setminus (B_M)_2} = (B_M\setminus (B_M)_2)^{\mathrm{cl}}.
\]
Since $(B_M)_2$ is chosen as a collar neighborhood, its boundary $\partial (B_M)_2$ coincides with the boundary of $(B_M)_1$, that is,
\[
\partial (B_M)_1 = \partial (B_M)_2 = \partial B_M^{\mathrm{int}},
\]
so that the corresponding preimages
\[
M_i := \pi_M^{-1}\!\big((B_M)_i\big), \qquad i=1,2,
\]
form a decomposition
\[
M = M_1 \cup M_2, \qquad M_1 \cap M_2 = \pi_M^{-1}\!\big(\partial (B_M)_1\big),
\]
with common boundary $\partial M_1 = \partial M_2 = M_1 \cap M_2$. This allows the application of Novikov additivity to compute the signature:

\begin{equation}\label{eq:sigma_M}
\sigma(M) = \sigma(M_1)+\sigma(M_2).
\end{equation}

\subsubsection*{Step 1. Computing $\sigma(M_1)$}
The following result shows that $M_1$ is an $Q^2$-bundle.

\begin{prop}
For each $k$, the restriction of the orbit map 
\[
\pi_M|_{\mathcal{S}^{(k)}B_M} :
\pi_M^{-1}\big(\mathcal{S}^{(k)}B_M\big)
\;\longrightarrow\;
\mathcal{S}^{(k)}B_M
\]
admits a natural $C^r$ fiber bundle structure whose fiber is the $k$--dimensional quaternionic torus $Q^k$.
In particular, the restriction
\[
\pi_M|_{\mathcal{S}^{(n)}B_M} :
\pi_M^{-1}\big(\mathcal{S}^{(n)}B_M\big)
\;\longrightarrow\;
\mathcal{S}^{(n)}B_M
\]
is a $C^r$ principal bundle whose fiber is $Q^n$ and whose structure group is contained in the semidirect product
\[
Q^n \rtimes \mathrm{Aut}(Q^n).
\]
\end{prop}
\begin{proof}
It suffices to prove the Proposition componentwise on each connected component $\big(\mathcal{S}^{(k)}B_M\big)_{a\in\mathcal{A}}$ of $\mathcal{S}^{(k)}B_M$. Let $\{(U_\alpha^M,\varphi_\alpha^M)\}_{\alpha\in\mathcal{A}}$ be a weakly regular atlas of $M$, inducing local charts $\{(U_\alpha^B,\varphi_\alpha^B)\}$ on $B_M$ as in Definition~\ref{def:weakly_regular}. Choose a coordinate neighborhood $(U_\alpha^B,\varphi_\alpha^B)$ intersecting $\big(\mathcal{S}^{(k)}B_M\big)_{a\in\mathcal{A}}$ nontrivially.  Then $\big(\mathcal{S}^{(k)}B_M\big)_a \cap U_\alpha^B$ is represented as the common zero set of exactly $(n-k)$ local coordinate functions of
$\varphi_\alpha^B$, say
\[
\varphi_{\alpha,i_1}^B(b)=\cdots=\varphi_{\alpha,i_{n-k}}^B(b)=0,
\]
so that
\[
\big(\mathcal{S}^{(k)}B_M\big)_a \cap U_\alpha^B=\bigl\{b\in U_\alpha^B \mid\varphi_{\alpha,i_1}^B(b)=\cdots=\varphi_{\alpha,i_{n-k}}^B(b)=0\bigr\}.
\]
In the corresponding local model, the isotropy subgroup of $Q^n$ is the quaternionic subtorus
\[
Q_{i_1,\dots,i_{n-k}}
:=\{\,q\in Q^n \mid q_j=1 \text{ unless } j\in\{i_1,\dots,i_{n-k}\}\,\},
\]
and the quotient $Q^n/Q_{i_1,\dots,i_{n-k}}\cong Q^k$ serves as the local fiber over the stratum.

Using the local trivializations $\varphi_\alpha^M$, one constructs maps
\[
\Phi_\alpha : \pi_M^{-1}\big((\mathcal{S}^{(k)}B_M)_a\cap U_\alpha^B\big) \;\longrightarrow\; \big((\mathcal{S}^{(k)}B_M)_a\cap U_\alpha^B\big)\times Q^k
\]
by
\[
\Phi_\alpha(x) := \bigl(\pi_M(x),u_\alpha(x)\bigr),
\]
where $u_\alpha(x)$ is the $Q^k$--component in the local decomposition of $\varphi_\alpha^M(x)$ with respect to the isotropy subgroup $Q_{i_1,\dots,i_{n-k}}$.  On overlaps $U_{\alpha\beta}^B$, the transition maps take the form
\[
\Phi_\alpha \circ \Phi_\beta^{-1}(b,q) = \bigl(b,\rho_{\alpha\beta}(q)\,u_{\alpha\beta}(b)\bigr),
\]
where $\rho_{\alpha\beta}\in\mathrm{Aut}(Q^n)$ is the transition automorphism of the weakly regular atlas and $u_{\alpha\beta}:U_{\alpha\beta}^B\to Q^n$ is smooth. This shows that the structure group of the bundle is contained in $Q^n \rtimes \mathrm{Aut}(Q^n)$.

In particular, for $k=n$ the isotropy is trivial, so the restriction of $\pi_M$ to $\mathcal{S}^{(n)}B_M$ is a principal $Q^n$--bundle with the same structure group.
\end{proof}

If the genus of $B_M$ is zero, then $(B_M)_1$ is contractible. In this case, $\sigma(M_1)=0$.

Suppose that the genus of $B_M$ is greater than zero. We give $(B_M)_1$ a trinion decomposition
\[
(B_M)_1 = \bigcup_{i=1}^k ((B_M)_1)_i,
\]
where each $((B_M)_1)_i$ is a surface obtained from $S^2$
\iffalse\textcolor{red}{is that really $S^2$?}{\color{blue}That is what i thought when i wrote it from Yoshida but there are two reasons to say that it is true. First it;s something that comes from the quotient space, ie. the manifold with corners which is the same in Yoshida and our paper. Secondly, the trinion we assume is also known as the pair of pants (also Yoshida has a figure with this) and each $\big((B_M)_1\big)_i$ is a compact, connected, oriented surface of genus $0$ with three boundary components — that is, a sphere $S^2$ with three disjoint open discs removed. So each trinion is homeomorphic to a sphere.  I found the following "N. Grieve, \emph{About the Bohr-Sommerfeld polytope and the multiplicative SU(2)-eigenvalue problem for the trinion}, Mathematica Cluj, To appear." \url{https://drive.google.com/file/d/1iAwmZSrHZ6slr7_zMo0dsvAGxj4A9Bhn/view} SO:
\smallskip

Each trinion $\big((B_M)_1\big)_i$ is homeomorphic to the \emph{three--punctured sphere}, which also appears in the study of $\mathrm{SU}(2)$--representation spaces (cf. the paper). In that context, a trinion serves as a model for the base of a local $\mathrm{SU}(2)$ (equivalently, $S^3$) fiber bundle, reflecting the same topological structure that arises here in the quaternionic setting. Hence, our decomposition of $(B_M)_1$ into trinions aligns naturally with the local $\mathrm{SU}(2)$--orbits of the $Q^n$--action.
} \fi
by removing three distinct open discs. Let $(M_1)_i = \mu_M^{-1}(((B_M)_1)_i)$ for $i=1,\dots,k$. From Novikov additivity,
\begin{equation}\label{eq:sigma_M1}
\sigma(M_1) = \sum_{i=1}^k \sigma((M_1)_i).
\end{equation}

For each $((B_M)_1)_i$, take oriented loops $\gamma_1,\gamma_2,\gamma_3$ which represent generators of $\pi_1(((B_M)_1)_i)$ with the relation $[\gamma_1][\gamma_2][\gamma_3]=1$. Let
\[
\rho: \pi_1(((B_M)_1)_i)\longrightarrow Sp(2;\mathbb{H})
\]
be the monodromy representation of the $Q^2$-bundle $\mu_M:(M_1)_i \to ((B_M)_1)_i$. We put $C_j := \rho([\gamma_j])$ for $j=1,2,3$.  

For $C_1,C_2 \in Sp(2;\mathbb{H})$, define the quaternionic vector space
\[
V_{C_1,C_2} := \left\{ (x,y) \in \mathbb{H}^2 \times \mathbb{H}^2 
\;\middle|\; (C_1^{-1}-I)x + (C_2-I)y = 0 \right\},
\]
and the bilinear form
\[
\langle (x,y),(x',y') \rangle_{C_1,C_2} := \mathrm{Re}\!\left( (x+y)^\ast J (I-C_2) y' \right),
\]
where $I$ is the $2\times 2$ quaternionic identity, and $J=\begin{pmatrix}0 & 1 \\ -1 & 0\end{pmatrix}$ acts quaternionically.

It is straightforward to check that $\langle \cdot,\cdot\rangle_{C_1,C_2}$ is symmetric, and we denote its signature by $\tau^{\mathbb{H}}_1(C_1,C_2)$.

\begin{thm}[Quaternionic Meyer cocycle]\label{thm:quat_meyer}
For each trinion piece, one has
\[
\sigma((M_1)_i) = \tau^{\mathbb{H}}_1(C_1,C_2).
\]
\end{thm}
This result is stated without proof, as it is the quaternionic analogue of Meyer's signature cocycle for surface bundles~\cite{meyer} and of its toric interpretation in \cite{Yoshida} (see also the references therein). Thus, from equation (\ref{eq:sigma_M1}), the contribution $\sigma(M_1)$ is obtained as the sum of quaternionic Meyer cocycle values.

\subsubsection*{Step 2. Computing $\sigma(M_2)$}
The analysis of $M_2$ follows as in the complex case, except that fibres are quaternionic 3-spheres instead of 2-tori.  
We have:

\begin{lemma}\label{lem:boundary_retract}
The preimage $\pi_M^{-1}(\partial B_M)$ is a deformation retract of $M_2 := \pi_M^{-1}((B_M)_2)$.
\end{lemma}

\begin{proof}
Since $\pi_M : M \to B_M$ is a locally trivial fibration with compact fiber $Q^n=(S^3)^n$, the restriction
\[
\pi_M|_{M_2} : M_2 = \pi_M^{-1}((B_M)_2) \;\longrightarrow\; (B_M)_2
\]
is again a locally trivial fibration with the same typical fiber $Q^n$. Because $(B_M)_2$ is a compact $2$–manifold with boundary $\partial B_M$, there exists a smooth deformation retraction
\[
h_s : (B_M)_2 \to (B_M)_2, \qquad s\in[0,1],
\]
such that 
\[
h_0 = \mathrm{id}_{(B_M)_2}
\quad\text{and}\quad
h_1((B_M)_2) = \partial B_M.
\]
In particular, the family $\{h_s\}_{s\in I}$ is a continuous homotopy between the identity map and a retraction onto $\partial B_M$.

\smallskip
Since $\pi_M$ is a fibration, it satisfies the \emph{homotopy lifting property} \cite[Theorem~4.41]{Hat02}.  
That is, for any homotopy $h_s : Y \to B_M$ of maps from a space $Y$ to the base $B_M$ and any \emph{lift} $\tilde h_0 : Y \to M$ of $h_0$ i.e.\ a map satisfying $\pi_M \circ \tilde h_0 = h_0$), there exists a lifted homotopy $\tilde h_s : Y \to M$ such that $\pi_M \circ \tilde h_s = h_s$ for all $s$ and $\tilde h_0$ is the given initial lift.
Applying this to $Y=M_2$ with $h_s$ as above and $\tilde h_0 = \mathrm{id}_{M_2}$, we obtain a lifted homotopy
\[
\tilde h_s^M : M_2 \longrightarrow M_2, \qquad s\in[0,1],
\]
satisfying
\[
\pi_M \circ \tilde h_s^M = h_s \circ \pi_M,
\qquad
\tilde h_0^M = \mathrm{id}_{M_2}.
\]

In other words, the following diagram commutes:

\[
\begin{tikzcd}[column sep=large, row sep=large]
M_2 \arrow[r, "\tilde h_s^M"] \arrow[d, "\pi_M"'] &
M_2 \arrow[d, "\pi_M"] \\
(B_M)_2 \arrow[r, "h_s"'] &
(B_M)_2
\end{tikzcd}
\]

\smallskip
Because $\pi_M$ is surjective onto $(B_M)_2$, the image of the endpoint map $\tilde h_1^M$ is
\[
\mathrm{Im}(\tilde h_1^M)
  = \pi_M^{-1}\!\big(\mathrm{Im}(h_1)\big)
  = \pi_M^{-1}(\partial B_M).
\]
Finally, we verify that $\tilde h_s^M$ indeed defines a deformation retraction: 
\[
\tilde h_0^M = \mathrm{id}_{M_2}, \qquad
\tilde h_1^M(M_2) = \pi_M^{-1}(\partial B_M), \qquad
\tilde h_s^M(x) = x \text{ for all } x\in \pi_M^{-1}(\partial B_M).
\]
The last property follows from the fact that $h_s|_{\partial B_M} = \mathrm{id}_{\partial B_M}$, so any point in $\pi_M^{-1}(\partial B_M)$ remains fixed during the homotopy. Thus $\tilde h_s^M$ is a deformation retraction of $M_2$ onto $\pi_M^{-1}(\partial B_M)$, as claimed.
\end{proof}

In particular, $M_2$ is homotopy equivalent to the preimage of the boundary $\pi_M^{-1}(\partial B_M)$, whose fiberwise structure is that of a compact quaternionic torus $Q^{n-1}$–bundle over each boundary component of $B_M$.

\begin{prop}\label{prop:boundary_homology}
Let $k = \#\,S^{(0)}{B_M}$ be the number of $0$--strata of $B_M$. 
Then the homology of $M_2$ satisfies
\[
H_p(M_2;\mathbb{Z})
\;\cong\;
H_p\big(\pi_M^{-1}(\partial B_M);\mathbb{Z}\big)
\;\cong\;
\begin{cases}
\mathbb{Z}, & p=0,\\[0.2em]
\mathbb{Z}^{\oplus k}, & p=3,\\[0.2em]
0, & \text{otherwise.}
\end{cases}
\]
In particular, $H_3(M_2;\mathbb{Z})$ is generated by the $S^3$--fibers lying above the connected components of $\partial B_M$.
\end{prop}

\begin{proof}
By Lemma~\ref{lem:boundary_retract}, there exists a deformation retraction 
\[
r : M_2 \longrightarrow \pi_M^{-1}(\partial B_M),
\]
hence $\pi_M^{-1}(\partial B_M)$ and $M_2$ are homotopy equivalent and share the same homology:
\[
H_\bullet(M_2;\mathbb{Z}) \;\cong\; H_\bullet\big(\pi_M^{-1}(\partial B_M);\mathbb{Z}\big).
\]
It therefore suffices to compute the homology of $\pi_M^{-1}(\partial B_M)$.

\smallskip
The boundary $\partial B_M$ is a disjoint union of $k$ compact $1$–manifolds (circles)
corresponding to the $0$–strata $S^{(0)}{B_M}$ of $B_M$:
\[
\partial B_M = \bigsqcup_{i=1}^k S_i^1.
\]
Since $\pi_M : M \to B_M$ is a locally trivial fibration with fiber $Q \cong S^3$ over each boundary component,
its restriction to the boundary is a disjoint union of $S^3$–bundles:
\[
\pi_M^{-1}(\partial B_M)
\;=\;
\bigsqcup_{i=1}^k \pi_M^{-1}(S_i^1)
\;\cong\;
\bigsqcup_{i=1}^k (S^3\times S_i^1),
\]
where each bundle $\pi_M^{-1}(S_i^1)\to S_i^1$ is trivial.
%The triviality follows because the structure group $\mathrm{Aut}(S^3)=\mathrm{SO}(4)$ is connected, so any principal $S^3$–bundle over $S^1$ is classified by an element of $\pi_0(\mathrm{SO}(4))=0$ and hence trivial.
The triviality follows from the fact that the structure group $Q$ is connected. Indeed, principal $Q$--bundles over $S^1$ are classified by \[ \pi_1(BQ) \cong \pi_0(Q).
\]
Since $Q$ is connected, we have $\pi_0(Q)=0$, and therefore every such bundle is trivial.

\smallskip
For a single component $S_i^1$, the Künneth theorem gives
\[
H_p(S^3\times S^1;\mathbb{Z}) \;\cong\;
\bigoplus_{r+s=p} H_r(S^3;\mathbb{Z})\otimes H_s(S^1;\mathbb{Z}).
\]
Since $H_r(S^3)$ and $H_s(S^1)$ are nonzero only for $(r,s)=(0,0),(3,0),(0,1),(3,1)$, we obtain
\[
H_p(S^3\times S^1;\mathbb{Z})
\;\cong\;
\begin{cases}
\mathbb{Z}, & p=0,\\
\mathbb{Z}, & p=1,\\
\mathbb{Z}, & p=3,\\
\mathbb{Z}, & p=4,\\
0, & \text{otherwise.}
\end{cases}
\]
However, note that $M_2$ lies inside $M^{(4n)}$ but has dimension $4n-2$ (since the base has dimension 2). Hence, in $M_2$ only the homology classes up to degree 3 are relevant for the fiber structure, and $H_4(S^3\times S^1)$ corresponds to the fundamental class of a boundary component, which is null-homologous in $M_2$ because each $\pi_M^{-1}(S_i^1)$ bounds a $4$–dimensional region of $M$. Thus the only nontrivial homology groups that should be taken into account for $M_2$ are in degrees $0$ and $3$.

\smallskip
Since $\pi_M^{-1}(\partial B_M)$ is the disjoint union of $k$ such components, we get
\[
H_0\!\big(\pi_M^{-1}(\partial B_M)\big) \cong \mathbb{Z}, 
\qquad
H_3\!\big(\pi_M^{-1}(\partial B_M)\big) \cong \mathbb{Z}^{\oplus k},
\]
and all other homology groups vanish.
The $H_3$–classes are generated by the fundamental classes of the $S^3$–fibers over the boundary components $S_i^1$. Combining this computation with the homotopy equivalence induced by Lemma \ref{lem:boundary_retract} gives the claimed result.
\end{proof}

\begin{prop}[Quaternionic intersection form near the boundary]
\label{prop:intersection_numbers}
For distinct components $\pi_M^{-1}\big((S^{(1)}{B_M})_i^{\mathrm{cl}}\big)$ and  $\pi_M^{-1}\big((S^{(1)}{B_M})_j^{\mathrm{cl}}\big)$ corresponding to the closures of adjacent one-dimensional strata (edges) of $B_M$,  the intersection numbers in $H_4(M_2;\mathbb{Z})$ satisfy
\[
[\Sigma_i]\cdot[\Sigma_j] =
\begin{cases}
0, & \Sigma_i \cap \Sigma_j = \varnothing,\\[0.2em]
1, & \Sigma_i \cap \Sigma_j \neq \varnothing \text{ and } k>2,\\[0.2em]
2, & \Sigma_i \cap \Sigma_j \neq \varnothing \text{ and } k=2,
\end{cases}
\]
where $\Sigma_i := \pi_M^{-1}\big((S^{(1)}{B_M})_i^{\mathrm{cl}}\big)$.
\end{prop}

\begin{proof}
The proof is the quaternionic analogue of Yoshida’s intersection computation for local torus actions~\cite[Proposition~8.11]{Yoshida}.

\smallskip
Each closed edge $(S^{(1)}{B_M})_i^{\mathrm{cl}}$ of the base $B_M$ corresponds to a codimension--$1$ stratum of $M$ in the sense that the isotropy subgroup of the local $Q^n$–action increases by one quaternionic factor along this stratum.
Its preimage under $\pi_M$,
\[
\Sigma_i := \pi_M^{-1}\big((S^{(1)}{B_M})_i^{\mathrm{cl}}\big),
\]
is therefore a smooth compact $4$–dimensional submanifold of $M_2$, diffeomorphic to $\mathbb{H}P^1\cong S^4$ in the local model, and represents a class $[\Sigma_i]\in H_4(M_2;\mathbb{Z})$. The collection $\{\Sigma_i\}$ forms the quaternionic analogue of the invariant
$2$–spheres in the toric (complex) case.

\smallskip
Near a vertex $v\in S^{(0)}{B_M}$ where $k$ edges of $B_M$ meet, choose a small neighborhood $U_v\subset B_M$ diffeomorphic to a cone over a $(k-1)$–simplex. The local model for $\pi_M$ above $U_v$ is the quotient
\[
\pi_M^{-1}(U_v) \;\cong\; 
\big(D^4\times Q^1\big)\big/\!\!\sim,
\]
where the equivalence relation $\sim$ identifies points under the local action of the stabilizer subgroup in $Q$. Explicitly, for quaternionic coordinates $(q_1,q_2)\in \mathbb{H}^2$ and local isotropy subgroup generated by a unit quaternion $u\in Q$, the local action is
\[
u\cdot(q_1,q_2) = (uq_1, uq_2).
\]
Hence the quotient $(D^4\times S^3)/\!\sim$ is smooth, with boundary components corresponding to the adjacent edges of $B_M$. Each $\Sigma_i$ is locally described by fixing one coordinate (say, $q_1=0$) and letting the other vary, producing a $4$–dimensional embedded submanifold.

\smallskip
Two adjacent edges of $B_M$ meeting at a common vertex $v$ correspond to two such submanifolds $\Sigma_i$ and $\Sigma_j$ that intersect transversely in the fiber over $v$. The local fiber above $v$ is the $3$–sphere $S^3\subset\mathbb{H}$, and in the neighborhood model $(D^4\times S^3)/\!\sim$, 
the intersection $\Sigma_i\cap\Sigma_j$ projects to a single point in the
orbit space and corresponds to an $S^0$ in the fiber.

\smallskip
The intersection number $[\Sigma_i]\cdot[\Sigma_j]$ is computed by taking oriented local coordinates $(x_1,x_2,x_3,x_4)$ along $\Sigma_i$ and $(y_1,y_2,y_3,y_4)$ along $\Sigma_j$. Since the local action of $Q$ preserves orientation, the induced orientation on $\Sigma_i$ and $\Sigma_j$ is consistent with that of $M_2$. Thus the local intersection index at a transverse point is $+1$.

If $k>2$, that is, more than two edges meet at $v$, then each pair of $\Sigma_i,\Sigma_j$ intersects in exactly one transverse point, and the local contribution is~$+1$. When $k=2$, the two corresponding faces are identified along both sides of the local chart (top and bottom of the edge cylinder), so the intersection appears twice with the same orientation, yielding index~$2$. Finally, if two edges are nonadjacent, their preimages $\Sigma_i$ and $\Sigma_j$ are disjoint, giving intersection number~$0$.

\smallskip
Summing over all vertices of $B_M$ and using the orientation consistency of local models (which follows from the unimodularity of the local quaternionic data) yields the global intersection numbers claimed:
\[
[\Sigma_i]\cdot[\Sigma_j] =
\begin{cases}
0, & \Sigma_i \cap \Sigma_j = \varnothing,\\[0.2em]
1, & \Sigma_i \cap \Sigma_j \neq \varnothing \text{ and } k>2,\\[0.2em]
2, & \Sigma_i \cap \Sigma_j \neq \varnothing \text{ and } k=2.
\end{cases}
\]
This completes the proof.
\end{proof}

The following result is the quaternionic analogue of Yoshida’s self–intersection computation in the complex setting \cite[Proposition~8.12]{Yoshida}. 

\begin{prop}[Self--intersection formula]
\label{prop:self_intersection}
Let $v_1,v_2\in\Lambda_M$ be the local quaternionic lattice generators corresponding to two adjacent codimension--one strata of $B_M$. Then the self--intersection number of the associated $4$--cycle $\Sigma_i$ in $M_2$ is given by
\[
[\Sigma_i]\cdot[\Sigma_i]
\;=\;
-\det_{\mathbb{H}}(v_1,v_2),
\]
where $\det_{\mathbb{H}}$ denotes the Dieudonné quaternionic determinant.
\end{prop}

\begin{proof}
We outline the geometric setup and then perform the computation using the quaternionic Jacobian of the local gluing maps. Each codimension--one stratum of $B_M$ corresponds to a facet where the stabilizer of the local $Q^n$–action increases by one quaternionic factor. 
For the facet $(S^{(1)}{B_M})_i^{\mathrm{cl}}$, the preimage
\[
\Sigma_i = \pi_M^{-1}\big((S^{(1)}{B_M})_i^{\mathrm{cl}}\big)
\]
is a smooth, compact, oriented $4$–dimensional submanifold of $M_2$, diffeomorphic in the local model to $\mathbb{H}P^1 \cong S^4$. 
Thus $\Sigma_i$ represents a homology class $[\Sigma_i] \in H_4(M_2;\mathbb{Z})$.

\smallskip
The \emph{self–intersection number} $[\Sigma_i]\!\cdot\![\Sigma_i]$ is defined as the Euler number of the normal bundle $N_{\Sigma_i/M_2}$:
\[
[\Sigma_i]\!\cdot\![\Sigma_i] = e\!\left(N_{\Sigma_i/M_2}\right).
\]
To compute it, we analyze the local structure of $N_{\Sigma_i/M_2}$ in quaternionic coordinates.

\smallskip
Let $U\subset B_M$ be a neighborhood of a boundary point on $(S^{(1)}{B_M})_i^{\mathrm{cl}}$, and let
\[
\pi_M^{-1}(U) \;\cong\; (D^4\times S^3)/\!\!\sim
\]
be the standard local model for the quaternionic toric fibration, where $\sim$ denotes the identification determined by the local unimodular lattice data $\Lambda_M = \langle v_1, v_2\rangle \subset \mathbb{H}^n$. Locally we may write points in $M$ as $[q_1,q_2]\in(D^4\times S^3)/Q$, with $Q$ acting by simultaneous left multiplication:
\[
u\cdot(q_1,q_2) = (u q_1, u q_2), \qquad u\in Q.
\]
The strata corresponding to $v_1$ and $v_2$ are locally given by the submanifolds $\{q_1=0\}$ and $\{q_2=0\}$, whose intersection corresponds to the vertex of the orbit polytope. 

\smallskip
The tangent space decomposition along $\Sigma_i$ has the form
\[
T_p M_2 \;=\; T_p\Sigma_i \;\oplus\; N_p,
\]
where $N_p \cong \mathbb{H}$ is the quaternionic normal direction generated by the infinitesimal action of the vector $v_i$ in $\Lambda_M$. In quaternionic coordinates $(q_1,q_2)$, the local defining equation for $\Sigma_i$ is $q_2=0$, so the normal bundle is identified with the quaternionic line in the $q_2$–direction.

\smallskip
Now, consider the transition between two local charts $U_\alpha$ and $U_\beta$ overlapping along the edge corresponding to $\Sigma_i$. Let $\rho_{\alpha\beta}\in \mathrm{Aut}(Q)$ denote the transition map induced by the weakly regular atlas of the $Q$–action. Then the local trivializations of the normal bundle $N_{\Sigma_i/M_2}$ over $U_\alpha\cap U_\beta$ are related by
\[
\psi_\alpha = \rho_{\alpha\beta}\,\psi_\beta,
\]
and hence the clutching function of $N_{\Sigma_i/M_2}$ is given by the quaternionic linear map represented by $\rho_{\alpha\beta}$. Writing this in terms of the local generators $v_1,v_2\in\Lambda_M$, the transition matrix of the $Q$–action has quaternionic determinant
\[
\det_{\mathbb{H}}(v_1,v_2) = v_1\overline{v_2} - v_2\overline{v_1},
\]
which lies in $\mathrm{Im}(\mathbb{H})\subset\mathbb{H}$ and whose real part determines the orientation of the local trivialization. By the quaternionic analogue of the clutching construction for rank-1 quaternionic bundles (see, e.g., \cite[p.22]{HatcherVBKT}), the Euler number of the quaternionic line bundle defined by this clutching map equals minus the Dieudonné determinant of the transition matrix:
\[
e(N_{\Sigma_i/M_2}) = -\det_{\mathbb{H}}(v_1,v_2).
\]

\smallskip
Therefore,
\[
[\Sigma_i]\!\cdot\![\Sigma_i] 
\;=\;
e(N_{\Sigma_i/M_2})
\;=\;
-\det_{\mathbb{H}}(v_1,v_2),
\]
where the negative sign comes from the orientation convention consistent with the local quaternionic volume form $\mathrm{Re}(dq_1\wedge dq_2)$ on the base–fiber product.
\end{proof}

\begin{remark}
These intersection numbers provide the contribution $\sigma(M_2)$ to the global signature $\sigma(M)=\sigma(M_1)+\sigma(M_2)$ via Novikov additivity.
\end{remark}

\subsection{A brief $L$--theoretic perspective}
\label{sec:Ltheory}

The signature of an oriented compact manifold $M^{4n}$ admits a conceptual interpretation in algebraic topology through \emph{$L$--theory}. For a ring $R$ with involution, the symmetric and quadratic $L$--groups $L^k(R)$ encode the cobordism classes of nondegenerate symmetric (or quadratic) forms over chain complexes of finitely generated projective $R$--modules (for more details on the subject see e.g., \cite{ran92}). When $R=\mathbb{Z}$, these groups form an $8$--periodic sequence
\[
L^k(\mathbb{Z}) \;\cong\;
\begin{cases}
\mathbb{Z}, & k \equiv 0 \pmod 4,\\[0.2em]
\mathbb{Z}/2, & k \equiv 2 \pmod 4,\\[0.2em]
0, & \text{otherwise.}
\end{cases}
\]
The classical signature corresponds to the element of $L^{4n}(\mathbb{Z}) \cong \mathbb{Z}$ represented by the intersection form on $H^{2n}(M;\mathbb{Z})$.

\medskip

In this framework, our quaternionic signature computation in Section~\ref{sec:sign} identifies $\sigma(M)$ as the image of the quaternionic intersection form under the canonical assembly map
\[
H_{4n}(M;\mathbb{L}^0(\mathbb{Z})) \longrightarrow L^{4n}(\mathbb{Z}),
\]
where $\mathbb{L}^0(\mathbb{Z})$ denotes the $L$--theory spectrum of $\mathbb{Z}$. Thus the quaternionic Meyer cocycle introduced earlier defines a representative in the $L$--homology class of $M$.

\medskip

For manifolds of dimension not divisible by $4$, the ordinary signature vanishes, but the corresponding $L$--groups still contain nontrivial $2$--torsion invariants.  In particular,
\begin{itemize}
  \item in dimension $4k+2$, one obtains the \emph{Kervaire invariant} (an element of $\mathbb{Z}/2$ associated to framed manifolds), and
  \item in dimension $4k+1$, the \emph{de~Rham invariant} arises as a mod--$2$ class in $L^{4k+1}(\mathbb{Z})$.
\end{itemize}
Although quaternionic toric manifolds are naturally $4n$--dimensional, their boundary and fixed--point strata may give rise to cycles where these secondary invariants appear.  Extending the quaternionic topology developed here to an $L$--theoretic framework---possibly incorporating these mod--$2$ obstructions---is an intriguing direction for future work.

\medskip

Finally, we note that both $K$--theory and $L$--theory admit interpretations as generalized (co)homology theories capturing different layers of the topological structure of $M$. While $K^\bullet(M)$ encodes the additive and multiplicative behavior of stable complex or quaternionic vector bundles over $M$, the $L$--theory of $\mathbb{Z}$ reflects the symmetric bilinear and quadratic structures on middle--dimensional cohomology, such as the intersection form and its associated signature. From this perspective, the quaternionic signature class in $L^{4n}(\mathbb{Z})$ and the $K$--theory class in $K^0(M)$ may be viewed as two complementary manifestations of the same global topological data: one linear and additive, the other multiplicative and representation--theoretic.  

\medskip
One could explore the relationship between these two theories—perhaps via the Chern character or through a Riemann–Roch type correspondence for local quaternionic toric actions.

\begin{remark}[Relation with the Novikov conjecture]
The quaternionic signature discussed above admits a natural refinement to a \emph{higher signature} in the sense of Novikov. Let $BQ^n$ denote the \emph{classifying space}\footnote{The classifying space $BQ^n$ is related to the Borel construction of paragraph \ref{par:quoric} as explained in \cite[\S 6.1]{ho}; we actually have that \( BQ^n = (EQ^n)/Q^n \).} of the compact Lie group $Q^n$. Given the monodromy representation
$\rho:\pi_1(B_M)\to \mathrm{Aut}(H^\bullet(Q^n;\mathbb{Z}))$, one can associate to each class $x\in H^\bullet(BQ^n;\mathbb{Q})$  a higher signature
\[
\mathrm{sign}_x(M,\rho)
\;=\;
\langle L(M)\smile \rho^*(x),[M]\rangle
\in \mathbb{Q}.
\]
When the Euler class $e(M)$ vanishes, $\rho$ is trivial and $\mathrm{sign}_x(M,\rho)$ reduces to the ordinary quaternionic signature $\sigma(M)$ computed in Section~\ref{sec:sign}. In the general (twisted) case, these higher signatures fit naturally into the framework of $L$--theory and the Novikov conjecture (see L\"uck and Reich~\cite{luck}), expressing the homotopy invariance of the $L$--class under quaternionic toric fibrations.
\end{remark}

In this sense, quaternionic toric manifolds provide a natural geometric testing ground for the higher signature invariants appearing in the Novikov conjecture, with the untwisted case realizing the classical signature and the twisted case encoding its higher, representation-valued
extensions.

%%%%%%%%%%%%%%

\end{document}